\documentclass[a4paper,11pt]{amsart}

\usepackage[utf8]{inputenc}  
\usepackage[T1]{fontenc}  
\usepackage{amscd,amssymb,amsmath,latexsym,graphicx,epsfig,color,url}
\usepackage[active]{srcltx}
\usepackage[all]{xy}
\usepackage{hyperref}

\usepackage{todonotes} 
\hypersetup{breaklinks=true}

\usepackage{todonotes}

\newcommand{\MV}{\operatorname{MV}}

\newcommand{\vol}{\operatorname{vol}}

\newcommand{\Hom}{\operatorname{Hom}}

\newcommand{\Res}{\operatorname{Res}}

\newcommand{\C}{\mathbb{C}}

\newcommand{\N}{\mathbb{N}}
\renewcommand{\P}{\mathbb{P}}
\newcommand{\Q}{\mathbb{Q}}
\newcommand{\R}{\mathbb{R}}
\newcommand{\T}{\mathbb{T}}
\newcommand{\Z}{\mathbb{Z}}

\newcommand{\cA}{{\mathcal A}}

\newcommand{\cD}{{\mathcal D}}

\newcommand{\cO}{{\mathcal O}}

\newcommand{\bfa}{{\boldsymbol{a}}}
\newcommand{\bfb}{{\boldsymbol{b}}}

\newcommand{\bff}{{\boldsymbol{f}}}

\newcommand{\bfm}{{\boldsymbol{m}}}

\newcommand{\bft}{{\boldsymbol{t}}}

\newcommand{\bfx}{{\boldsymbol{x}}}
\newcommand{\bfy}{{\boldsymbol{y}}}
\newcommand{\bfz}{{\boldsymbol{z}}}
\newcommand{\bfsigma}{{\boldsymbol{\sigma}}}

\newcommand{\bfF}{{\boldsymbol{F}}}

\newcommand{\bfxi}{{\boldsymbol{\xi}}}
\newcommand{\bfchi}{{\boldsymbol{\chi}}}

\newcounter{thm}
\numberwithin{thm}{section}
\numberwithin{equation}{section}

\theoremstyle{definition}
\newtheorem{definition}[thm]{Definition}

\newtheorem{remark}[thm]{Remark}
\newtheorem{example}[thm]{Example}

\newtheorem{question}[thm]{Question}

\theoremstyle{plain}
\newtheorem{lemma}[thm]{Lemma}
\newtheorem{proposition}[thm]{Proposition}
\newtheorem{theorem}[thm]{Theorem}
\newtheorem{corollary}[thm]{Corollary}

\begin{document}

\title[Toric Euler-Jacobi theorem]{Toric Euler-Jacobi vanishing theorem and zeros at infinity }

\author{Carlos D'Andrea}
\address{Departament de Matem\`atiques i Inform\`atica, Universitat de Barcelona. Gran
Via 585, 08007 Barcelona, Spain \&  Centre de Recerca Matem\`atica, Edifici C, Campus Bellaterra, 08193 Bellaterra, Spain
}
\email{cdandrea@ub.edu}
\urladdr{http://www.ub.edu/arcades/cdandrea.html}

\author{Alicia Dickenstein}
\address{Departamento de Matem\'atica, FCEN, Universidad de Buenos Aires, Ciudad Universitaria - Pab. I, (1428) Buenos Aires, Argentina}
\email{alidick@dm.uba.ar}
\urladdr{http://mate.dm.uba.ar/~alidick}

\begin{abstract}

Residues appear naturally in various questions in complex and algebraic geometry: interpolation, duality, representation problems, and obstructions. 
The first global vanishing result in the projective plane, known as the Euler-Jacobi theorem, was established by Jacobi in 1835. 
In the toric case, the input is a system of $n$ Laurent sparse polynomials with fixed Newton polytopes, and the first version 
of the Euler-Jacobi toric vanishing theorem for residues in the $n$-torus is due to Khovanskii in 1978, under restrictive genericity
 assumptions. In  this paper, we provide geometric conditions on the input Newton polytopes to ensure that this global vanishing is 
 equivalent to the existence of zeros at infinity in the associated compact toric variety. We relate these conditions to the dimension 
 at the toric critical degree of the quotient of the Cox ring by the ideal generated by the (multi)homogenizations of the input polynomials. 
 We also relate the existence of zeros at infinity to interpolation questions.

\end{abstract}


\maketitle

\section{Introduction}\label{sec:intro}

\subsection{A bit of history and our aim}\label{ssec:history}

Jacobi started his 1835 paper in latin {\em New algebraic theorems on systems of two equations in two unknowns}~\cite{Jacobi}  
with the following mention to {\bf Euler's theorem}: ``Of the theorems expressed in algebraic terms, there is hardly any more useful where equations are involved than the following, well-known one:  

\smallskip

If $X$ is a polynomial in $x$, we have that
\begin{equation} \label{eq:Euler} \sum \left( \frac U {\frac {dX} {dx}}\right) \, = \, 0,\end{equation}
if the sum runs over all roots $x$ of the equation $X=0$, and $U$ is a polynomial in $x$ of degree two less than the degree of $X$.''

\smallskip

His aim was to prove that this theorem can be extended to a system of two algebraic equations with two unknowns. We state his main result, translated into English.

\smallskip

{\bf Theorem (Jacobi):}
 Let  $f, \varphi$ be polynomials in two variables $x,y$; let $F$ be any other polynomial of {\em degree smaller than
the sum of the degrees of $f$ and $\varphi$ minus $2$}; it will be
\begin{equation}\label{eq:Jacobi}
\sum \frac F {f'(x) {\varphi}'(y) - f'(y) {\varphi}'(x) }\, = \, 0,\end{equation}
(where) the sum runs over all values $x,y$ that are common roots  of the equations $f=0, \varphi=0$.

\smallskip

Note that the denominator is a way of writing the {\em Jacobian}  $J_{f, \varphi}$ of $f$ and $\varphi$, that Jacobi introduced
later in 1841! He didn't write two hypotheses. The first one is that the jacobian  {\em does not vanish}  
at the common roots, and this was also an omitted hypothesis in the statement of
Euler's theorem, where it is assumed that all roots of the input polynomial are simple (that is, the derivative does not vanish at any of the roots). 

 Khovanskii suggested to us the following nice proof of Euler's theorem. 
 Let $f$ be a univariate polynomial with complex coefficients of degree $d$, with distinct roots $\xi_1, \dots, \xi_d$. Thus, $f= a_d \prod_{i=1}^d (x - \xi_i)$, with $a_d\neq 0$.
Let $L_i(x) \, = \, \frac {\prod_{ j \not= i} (x- \xi_j)}{\prod_{ j \not= i}
(\xi_i - \xi_j)}$ be the associated Lagrange
interpolating polynomial  (of degree $d-1$).
For any polynomial $h$ with ${\rm deg}(h) \leq d-1$,
$h(x) = \sum_{i=1}^d h(\xi_i) \, L_i(x)$.
So, the coefficient of $x^{d-1}$ in this sum
should be $0$ if  $ \deg(h)$ is at most $d-2$. But this coefficient is precisely
$$ {0} = \sum_{i=1}^d h(\xi_i) \frac 1 {\prod_{ j \not= i} (\xi_i - \xi_j)} = a_d
\,
 \sum_{i=1}^d \frac {h(\xi_i)} {f'(\xi_i)}.$$
This result can be easily generalized to the case of possibly repeated roots by noting that
$\frac {h(\xi_i)} {f'(\xi_i)}$ equals $ 2 \pi i$ times the local {\em residue} of the differential form 
$\Phi_f(h) =\frac h f \, dx$ at $\xi_i$.  Another proof may be given using the basic result that the sum of local residues
at all the zeros of $f$ in $\mathbb C$ equals minus the residue at infinity of the extension to complex projective space $\P^1$ of $\Phi_f(h)$
and it is easy to see that this residue vanishes whenever $\deg(h) < d-1$.

The second missing hypothesis in Jacobi's statement for the case of two variables is that $f$ and $\varphi$ 
have the B\'ezout number  $\deg(f)\cdot \deg(\varphi)$ of isolated common roots in the complex plane, or equivalently, that $f$ and $\varphi$  do {\em not} have any
common zeros at infinity in $\mathbb {P}^2$, as we will deduce from Theorem~\ref{smt}.  
Again, for each isolated common root $(x_i,y_i)$ of $f$ and $\varphi$ the summand $\frac{F(x_i,y_i)}{J_{f, \varphi} (x_i,y_i)}$ equals $(2 \pi i)^2$
times the local residue at $(x_i,y_i)$ of the differential form $\frac F {f \cdot \varphi} \, dx \wedge dy$, see for instance~\cite{CD05a,GH78}. 
 Equality~\eqref{eq:Jacobi} asserts that the sum of local residues
vanishes for any polynomial $F \in \C[x,y]$ with $\deg(F) < \deg(f) + \deg(\varphi) -2$. Indeed, Euler and Jacobi theorems give vanishing conditions 
for the sum of residues in the projective  compactifications of $\C$ and $\C^2$. Chapter~5 in \cite{GH78} shows how classical geometric
theorems in the projective plane are a consequence of this global vanishing.

There is a natural generalization of Jacobi's vanishing theorem to the multivariate case: given 
$n$ polynomials $f_1, \dots, f_n \in \C[x_1, \dots, x_n]$ with isolated zeros $\xi_1, \dots, \xi_m$ in $\C^n$ 
such that the sum of the corresponding intersection multiplicities equals the product of their degrees (that is, such that the hypersurfaces
$(F_i=0)$ defined by the homogenizations $F_1, \dots, F_n$ of the input polynomials do not have any intersection point at infinity in $\P^n= \P^n(\C)$, 
 i.e. at $(x_0=0)=\P^n\setminus\C^n$), if $h\in\C[x_1,\ldots, x_n]$ has degree smaller than the sum of the total degrees of the $f_i$'s  minus $n,$ then the sum of local residues of the form $\frac h {f_1 \cdots f_n} \, dx_1 \wedge \dots \wedge dx_n$  over
$\xi_1, \dots, \xi_m$, is equal to zero. See for instance Corollary $4$ in \cite[Chapter 5]{AGV85} or Theorem 7.10 in \cite{kun08}.

To extend this result to the sparse setting, the role of $\P^n$ is played by a toric variety $X_P$ which will be described below. The polynomials 
$f_1,\ldots, f_n$ can be {\em (multi)homogenized}  to $F_1,\ldots, F_n$ in  the \textit{homogeneous coordinate ring} $S$ of $X_P$~\cite{cox95},  
as it will be recalled in Section~\ref{sec:1}, so that it has sense to consider their zeroes in $X_P$. However, in the general toric case the irrelevant ideal is in general  
not generated by a regular sequence and the commutative algebra 
associated to the corresponding multigrading is quite complicated.

The input data are $n$ {\em lattice polytopes} $P_1, \dots, P_n \subset \R^n$ . We will always assume that their Minkowski sum
 $ P \, = \, P_1 + \dots + P_n$ is an $n$-dimensional lattice polytope and we will denote by $P^\circ$ its interior with respect to the Euclidean topology. We then consider Laurent polynomials 
 $f_1, \dots, f_n \in \C[\bft^{\pm1}]:=\C[t_1^{\pm1},\ldots, t_n^{\pm1}]$ with monomials $\bft^{\bfa} = t_1^{a_1}\dots t_n^{a_n}$ in $P_1, \dots, P_n$:
\begin{equation}\label{eq:efes}
f_i \, = \, \sum_{\bfa\in P_i \cap \Z^n} c_{i{\bfa}}\bft^{\bfa}, \ i=1,\ldots, n,
\end{equation}
and denote by $V_T(\bff) =V_{(\C^*)^n}(f_1,\ldots, f_n)$ the variety of their common zeros in the $n$-torus $T = (\C^*)^n$, with  $\C^*:=\C\setminus\{0\}$.
 B\'ezout theorem is extended to the Bernstein-Kouchnirenko-Khvovanskii bound. 
Bernstein Theorem (\cite{ber75}) states that the degree $\deg(V_T(\bff))$  of the zero dimensional part of $V_T(\bff)$ (the sum of the multiplicities of all its isolated zeroes over $T$) is bounded by 
the mixed volume $\MV(P_1,\ldots, P_n),$ and it is equal to this number when the polynomials have generic coefficients (cf. Section 7.5 in~\cite{CLO05}).
 There is a natural compactification 
of the $n$-torus $T$ given by the compact toric variety $X_P$ associated with the normal fan of $P$.
Then,  $X_P \setminus T = \cup_i D_i$, with one divisor $ D_i$ at infinity associated  with any facet of $P$.

We briefly recall the definition of local and global residues and introduce the main ingredients in the toric case. We expand these definitions in Section~\ref{sec:3}.
Given a sequence of $n$ Laurent polynomials in $n$ variables $f_1,\ldots, f_n\in \C[\bft^{\pm1}]$ 
such that the variety $V_T(\bff)$ is finite,  and any  $h\in \C[\bft^{\pm1}]$, the \emph{global residue} $\mbox{Res}^T_\bff(h)$ in $T$  of the differential form
$\omega_h:=\frac{h}{f_1\ldots f_n}\frac{dt_1}{t_1}\wedge\ldots\wedge\frac{dt_n}{t_n},
$
is defined as the sum of the local Grothendieck residues of $\omega_h$ at each of the points of $V_T(\bff):$
\begin{equation}\label{eq:globalr}
\mbox{Res}^T_\bff(h):=\sum_{\bfxi\in V_T(\bff)}\mbox{Res}_{\bfxi}(\omega_h).
\end{equation}
The proper analytical definition of $\mbox{Res}_{\bfxi}(\omega_h)$ is given in \eqref{reslocal}. 
A more informal characterization of this operator is the following:
the \emph{jacobian in the torus} is defined as 
\begin{equation}\label{afj}
J^T_\bff:=\det\big(t_i\frac{\partial f_j}{\partial t_i}\big)_{1\leq i, j\leq n}.
\end{equation}
For $\bfxi\in V_T(\bff),\, \mbox{Res}_{\bfxi}(\omega_h)$  is equal to the value of $\frac{h(\bfxi)}{J^T_\bff(\bfxi)}$ if $\bfxi$ is a simple root of $V$, or to $\lim_{\varepsilon\to 0^+}\sum_{\bfxi_\varepsilon} \frac{q(\bfxi_\varepsilon)}{J^T_{\bff_\varepsilon} (\bfxi_\varepsilon)}$ of a deformed system $f_{\varepsilon,i}\to f_i,\,i=1,\ldots, n$, the sum being over all $\bfxi_\varepsilon\in V_T(\bff_\varepsilon)$ which converge to $\bfxi$ when $\varepsilon\to0.$ This limit always exists, see for instance \cite[\S 5.11]{AGV85}. There are also purely algebraic definitions of the residue, we refer the reader to \cite{kun08,BCRS96} for more in this direction.
Global residues in the torus and their properties have been studied in \cite{kho78b, sop07, CD97, CDS98}. They have important connections and applications in algebraic geometry,  polynomial system solving and elimination theory, see for instance \cite{CM96, EM00, EM04, CD05a, sop13, SZ16, BSST23, BCCFMM24, SY24, vill24}. 
 
We have the following generalization of the {Euler-Jacobi theorem} to the {toric} setting from~\cite[Theorem 2]{kho78} under
genericity hypotheses, then generalized in~\cite[Corollary $5$]{CD97}: 

\smallskip

\begin{theorem} {{\rm (}\bf Toric Euler-Jacobi theorem{\rm )}}\label{koo} 
Assume $V_T(\bff)$ is finite. If the sum of the local multiplicities
 at all the points $\xi \in V_T(\bff)$ equals $ MV(P_1, \dots, P_n)$ (equivalently, the closures of the 
 hypersurfaces $(f_i=0)$ do not have any common point of intersection at the divisors at infinity of $X_P$), then for any
Laurent polynomial $h \in \C[\bft^{\pm1}]$ with monomials in the interior $P^\circ$ of $P$, the sum of local residues
in~\eqref{eq:globalr} 
is equal to 0:
\[\mbox{Res}^T_\bff(h) \, = \, 0.\]
\end{theorem}

 \medskip

In this paper, we explore the {\em necessity} of the hypothesis of not {\em having zeros}
 at the toric (and in particular, at the projective) infinity, that is: {\it is the (toric) Euler-Jacobi
vanishing equivalent to this hypothesis}? 
Our main result is Theorem \ref{smt} where in particular we prove a general converse of  Theorem \ref{koo}.

\subsection{The converse of the Euler-Jacobi Theorem and other results}\label{ssec:cej}

The next cautionary example shows that the converse to Theorem~\ref{koo} is not true in full generality.

\begin{example} \label{ex:notindec}
Set $n=3$. We consider the lattice polytopes $P_1, P_2, P_3$  given as the convex hulls of the following configurations of lattice points: 
$\{(0,0,0),\,(1,0,0),\,(0,1,0)\}$ , $\{(0,0,0),\,(0,1,0),\,(1,0,0),\,(2,0,0),\,(1,1,0),\,(0,2,0)\},$ and $\{(0,0,0),\,(0,0,1),\,(0,0,2)\}$. Denote by $\Delta_2$ the standard simplex of $\R^2,$ i.e. the convex hull of $(0,0),\,(1,0)$ and $(0,1).$ Then, we have that $P_1=\Delta_2 \times \{0\},\, P_2=2\Delta_2\times \{0\},$ and $P_3=\{(0,0)\} \times [0,2]$. Thus, their Minkowski sum $P = P_1 +P_2 + P_3$ equals the cartesian product $3 \Delta_2 \times [0,2]$ and has dimension $3$. It is easy to verify that $\MV(P_1,P_2, P_3)=4.$

 Consider the following polynomials with these supports:
\begin{equation}\label{eq}
\begin{array}{ccl}
f_1&=& 1+t_1+t_2\\
f_2&=& 2 - t_1 + t_1^2 + t_2 + 2 t_1 t_2 + t_2^2\\
f_3&=& 2-3t_3+t_3^2.
\end{array}
\end{equation}
The system~\eqref{eq} has only  $2$ simple solutions in $(\C^*)^3:\, (1,-2,1)$ and $(1,-2,2).$   The total number of zeros (of the associated divisors) in the complete toric variety $X_P$ (recalled in Section~\ref{ssec:vtp}) is finite; indeed, it is straightforward to see that they have $4$ simple zeros in $X_P$.
The only lattice  point lying in the interior of $P_1+P_2+P_3$ is $(1,1,1),$ and the jacobian $J_{f_1, f_2, f_3}$ is equal to $-4t_3+6,$ which has opposite values in the two common zeroes of the system in the torus. This implies straightforwardly that the global residue in the torus of $t_1t_2t_3$ is equal to zero, despite the fact that there are zeroes at infinity of the system. 
\end{example}

We need to introduce the following definition.
\begin{definition}\label{def:polytopes}
A sequence of polytopes $P_1,\ldots, P_n$ is said to be \textit{essential} if for any $k=1,\ldots, n,$ and any $J\subset\{1,\ldots, n\}$ of cardinality $k$, the dimension of $\sum_{j\in J} P_j$ is at least $k$. 
A sequence of polytopes $P_1,\ldots, P_n$ is said to be \textit{indecomposable} if it is essential and for any $k=1,\ldots, n-1,$ and any $J\subset\{1,\ldots, n\}$ of cardinality $k$,  either the dimension of $\sum_{j\in J} P_j$ is not equal to $k,$ or if it this happens, there are no lattice points in the relative interior of this Minkowski sum with respect to the Euclidean topology.
\end{definition}

 For instance, the sequence is essential and indecomposable if all the polytopes are $n$-dimensional.   The fact that a sequence of polytopes is essential is equivalent to the fact that their mixed volume is positive (that is, that generic polynomials with these supports have a nonempty set of common roots), see for example Theorem~4.13 in~\cite{E96}. In particular, their Minkowski sum has maximal dimension $n$.

 We now state our main result. 

\begin{theorem}\label{smt}
Let $f_1,\ldots, f_n\in\C[\bft^{\pm1}]$ be as in \eqref{eq:efes} with  $P_1,\ldots, P_n$ indecomposable. Assume that $V_T(\bff)$ is nonempty and that the intersection of the closures of the hypersurfaces $(f_i=0)$ for $i=1,\dots,n$, has dimension zero in  the complete toric variety $X_P$ associated with $P=P_1+\dots +P_n$ and it is contained in the simplicial part of $X_P$ (as in Definition~\ref{def:simplicial}).
Then, the following are equivalent:
\begin{enumerate}
\item[i)] $\deg(V_T(\bff))=\MV(P_1,\ldots, P_n)$.
\item[ii)] For any $ h_0\in\C[\bft^{\pm1}]$ with support contained in $P^\circ$, 
$\mbox{Res}^T_{\bff}(h_0) =0$.
\item[iii)] There is no Laurent polynomial  $p_J\in  \C[\bft^{\pm1}]$ with support contained in $P^\circ$ such that $J^T_\bff\equiv p_J$ modulo the ideal 
$\langle f_1,\ldots, f_n\rangle$ in $\C[\bft^{\pm1}]$. 
\end{enumerate}
\end{theorem}

Note that there is a stronger condition in the hypothesis of Theorem \ref{smt} on the polytopes $P_i,\,i=1,\ldots, n,$ because we are requiring the sequence to be indecomposable, as opposed to the classical Euler-Jacobi in the torus (Theorem \ref{koo}) where the condition was  only  that the polytope  $P = P_1 + \dots + P_n$ is $n$-th dimensional. This strong condition is necessary as Example~\ref{ex:notindec} shows that item ii) does not necessarily imply item i) when it is not satisfied. 

\smallskip

 In another but related direction, it is known that for a system of homogeneous polynomials $G_0,\ldots, G_n\in\C[x_0,\ldots, x_n],$ their standard jacobian belongs to the ideal $\langle G_0,\ldots, G_n\rangle$ if and only if the variety they define in projective space $\P^n$ is empty~\cite{spo89,vas92} (see Section~\ref{sec:open}). In the sparse context, there is a notion of \textit{toric jacobian}  $J_\bfF$ and a more general \textit{discrete jacobian} $\Delta_{\bfF,\bfsigma}$ (see Section~\ref{sec:critical}, and in particular~\eqref{JF} and~\eqref{delta}) for a system of (toric) homogeneous polynomials $F_0,\ldots, F_n$ in the homogeneous coordinate ring $S$ of $X_P$ with ample degrees.
 
 We now present our second main result:
 
 \begin{theorem}\label{fonlyif}
Let $X$ be a complete toric variety and $F_i\in S_{\alpha_i}$ for $i=0,\ldots, n,$ where each degree $\alpha_i$ is ample. Assume that $V_X(F_0,\ldots, F_n)\neq\emptyset$ and
that  for some $i\in\{0,\ldots, n\}$, the zero set $V_X(F_0,\ldots, F_{i-1}, F_{i+1},\ldots, F_n)$ is finite and contained in the simplicial part of $X$. Then,
\begin{enumerate}
\item If $\bfsigma$ is a flag of cones in $\Sigma,$  then $\Delta_{\bfF,\bfsigma}\in \langle F_0,\ldots, F_n\rangle.$
\item If in addition each of the $\alpha_i$'s is an integer multiple of a fixed ample Cartier degree  $\alpha,$ then  $J_\bfF\in \langle F_0,\ldots, F_n\rangle.$
\end{enumerate}
\end{theorem}

 It turns out that the dimension of certain graded pieces of the quotient of the homogeneous coordinate ring $S$   by the ideal generated by the homogenizations  $F_1, \dots, F_n$ of
the input Laurent polynomials $f_1, \dots, f_n$ is related to the  dimensions of  graded pieces of the quotient by the ideal generated by $F_1, \dots, F_n$ and another (well chosen) homogeneous polynomial $F_0$ such that the associated $(n+1)$ hypersurfaces have no common zeros. If the variety $X_P$ is simplicial, local and global residues associated to $F_0, \dots, F_n$  can also be defined in $X_P$, and they satisfy suitable vanishing conditions, which will be crucial to our problem, see Section \ref{groth}.   These tools, combined with the description of the structure of local points in a zero-dimensional complete intersection in $X_P$, will allow us to prove our main results.

This paper is organized as follows: in Section~\ref{sec:1} we recall the concept of  toric varieties and fix some notation. Critical degrees and toric jacobians are recalled in Section~\ref{sec:critical}. Residues are introduced in Section~\ref{sec:3}. We then review recent results by Bender and Telen (cf. \cite{BT22}) on the local structure of affine zero dimensional complete intersections in toric varieties in Section \ref{sec:4}, and relate them to residues in toric varieties.   Section \ref{sec:proofs} will use the previous tools to prove a codimension one theorem (Theorem \ref{oa}). This result will allow us to prove both Theorem \ref{smt}  and Theorem \ref{fonlyif}. We conclude the paper with some open questions and directions of further research in Section~\ref{sec:open}.
 
 \medskip
\noindent\underline{\bf Acknowledgements:}
We are grateful to Mat\'ias Bender and Simon Telen for explaining to us the results of their paper \cite{BT22}, to Greg Smith for pointing us  references to the literature, to Alvaro Liendo for discussions about affine simplicial toric varieties, to Antonio di Scala for suggesting to us to look at Jacobi's paper~\cite{Jacobi} in latin,  to David Cox for 
providing us  precise references to the book~\cite{CLS11} in Section~\ref{ssec:divpol} and suggesting the proof of Proposition~\ref{prop:1}.
Part of this work was done while the second author visited the University of Barcelona
and the first author visited the University of Buenos Aires. 
 We are grateful to the support and stimulating atmosphere provided by both institutions. These research trips were funded partially by the Spanish MICINN research projects  PID2019-104047GB-I00 and PID2023-147642NB-I00. In addition, C. D'Andrea was supported by the Spanish State Research Agency, through the Severo Ochoa and Mar\'ia de Maeztu Program for Centers and Units of Excellence in R\&D (CEX2020-001084-M), and the European  H2020-MSCA-ITN-2019 research project GRAPES. 
A. Dickenstein was partially supported by UBACYT 20020220200166BA and CONICET PIP 11220200100182CO, Argentina.

\bigskip
\section{Basics on Toric varieties}\label{sec:1}
We review in this section the main notions about toric varieties that we will need. All the terminology and notation we use is taken from \cite{CLS11}. As there, our base field here is $\C,$ but any other algebraically closed field of characteristic zero would also work. 

\subsection{Toric varieties from fans}\label{ssec:tff}
Let $X$ be a {\em normal complete  toric variety} of dimension $n$. So, there is an $n$-dimensional lattice $N$, and $X = X_\Sigma$ is determined by a complete (polyhedral) fan $\Sigma$  in $N_\R:=N\otimes\R$ as follows.
 The toric variety $X_\Sigma$ is complete when the support of the fan equals $N_\R$. $X_\Sigma$ is {\em simplicial} if the fan $\Sigma$ is {\em simplicial}, i.e. if  for every cone $\sigma\in\Sigma,$ the minimal generators of $\sigma$ are linearly independent in $N_\R$.  The fact that $X_\Sigma$ is smooth can be also read from $\Sigma$ as this is equivalent to the fact that the generators of every maximal cone in the fan are a $\Z$-basis of $N$.

\begin{definition} \label{def:simplicial}
Let $X_\Sigma$ be a normal toric variety with fan $\Sigma$. Let ${\Sigma}'$ be the subfan of $\Sigma$ consisting of all its simplicial cones. The simplicial toric variety $X_{{\Sigma}'} \subseteq X_\Sigma$ is called the simplicial part of $X_\Sigma$.
\end{definition}
 
 \begin{remark} \label{rem:normal} We will be dealing with normal toric varieties associated with polyhedral fans, even if we omit to mention the property of normality along the text.
 \end{remark}

Let $M$ be the dual lattice of $N$, and denote by $\eta_1,\ldots, \eta_{n+r}$ ($r\geq1$) the primitive generators in $N$ of the $1$-dimensional cones in $\Sigma$. We have then that  $r$ is the rank of the {\em Chow group} $A_{n-1}(X)$, which can be characterized as the cokernel of the map 
$M\to \Z^{n+r}$ sending $m$ to $(\langle m,\eta_1\rangle ,\ldots, \langle m,\eta_{n+r}\rangle)$.
Let $S:=\C[x_1,\ldots, x_{n+r}]$ be the {\em homogeneous coordinate ring} of $X$, as defined in \cite{cox95}. The variety $X$ is  the (compatible) union of the affine open sets $X_\sigma$, with $\sigma$ being an $n$-dimensional cone in $\Sigma$ defined as follows. Set $\bfx^{\widehat{\sigma}}:=\prod_{\eta_j\notin\sigma} x_j.$  Then, $\C[X_\sigma]$ is the $0$-th graded piece of $S_{\bfx^{\widehat{\sigma}}},$ the localization of $S$ at the monomial $\bfx^{\widehat{\sigma}}$ (cf. \cite[Lemma 2.2]{cox95}), and the grading given by $A_{n-1}(X)$ as explained below. $X$ has  the torus $N\otimes_\Z\C^*=\mbox{Hom}(M,\C^*) \cong (\C^*)^n$ as an open dense subvariety, with coordinate ring $\C[\bft^{\pm1}]$ and the coordinatewise multiplication map of the torus by itself is extended to an action of the torus on $X$.  The coordinates $(t_1, \dots, t_n)$ in the torus are related to those given by $x_1,\ldots, x_{n+r}$ via 
\begin{equation}\label{tx}
t_i=\prod_{j=1}^{n+r} x_j^{\langle\eta_j,e_i\rangle}, \ i=1,\ldots, n.
\end{equation}

Each variable $x_j$ of $S$ corresponds to a generator $\eta_j$ of a $1$-dimensional ray in the fan $\Sigma$ and hence to a torus-invariant irreducible divisor $D_j$ of $X.$ We will  denote with $[D_j]\in A_{n-1}(X)$ the class of $D_j$ in the Chow group, $j=1,\ldots, n+r.$ As in \cite{cox95} we grade $S$ in such a way that the monomial
$\prod_{j=1}^{n+r}x_j^{a_j}$ has degree $[\sum_{j=1}^{n+r} a_jD_j]\in A_{n-1}(X).$ 
Set 
\begin{equation} \label{eq:rho}
\rho_0:=[\sum_{j=0}^{n+r}D_j]\in A_{n-1}(X),
\end{equation}
the sum of the degrees of all the variables.

Denote $M_\R:= M\otimes\R$. Associated to a divisor $D$ with class 
\begin{equation}\label{alpa}
[D]= \alpha = \sum_{j=1}^{n+r} a_j [D_j],
\end{equation} we can define the polytope
\begin{equation}\label{pd}
P_D = P_\alpha :=\{m\in M_\R:\,\langle m,\eta_j\rangle\geq-a_j,\,1\leq j\leq n+r\}\subset M_\R.
\end{equation}
Note that there is in general more than one way of writing the class $[D]$ as a linear combination of the classes of the divisors $[D_j]$ as in \eqref{alpa}, so the polytope $P_D$ is well defined up to translations by elements in $M$.
We have that by denoting with $S_\alpha$ the graded piece of $S$ of degree $\alpha$, then
\begin{equation}\label{kor}
S_\alpha \cong H^0(X,\cO_X(D))=\oplus_{m\in P_D\cap M}\C\cdot \bfchi^{m},
\end{equation}
where $\bfchi^m\in \mbox{Hom}(N\otimes_\Z\C^*,\C^*)$ denotes the  character in the torus associated with $m$, and $\cO_X(D)$ denotes the standard sheaf on $X$ associated to  $D$ (\cite[Proposition 1.1]{cox95}). This implies that any polynomial $F\in S_\alpha$ of degree $\alpha$ can be written as
\begin{equation}\label{eFe}
F=\sum_{m\in P_D\cap M}c_m\prod_{j=1}^{n+r} x_j^{\langle m,\eta_j\rangle+a_j},
\end{equation}
with $c_m\in\C.$

Recall that a divisor $D=\sum_{j=1}^{n+r} a_j D_j$ is {\em Cartier} if and only if for all $\sigma\in\Sigma$ generated by $\{\eta_{i_1},\ldots, \eta_{i_k}\},$ there exists $m_\sigma\in M$ such that $\langle m_\sigma,\eta_{i_j}\rangle=-a_{i_j}$ for $j=1,\ldots, k$  (\cite[Theorem 4.2.8]{CLS11}). For such divisors, one can \textit{dehomogenize} homogeneous polynomials in $S$ to elements in  $\C[X_\sigma]$ for every maximal $\sigma\in\Sigma$ as follows:
\begin{equation}\label{deho}
F\in S_\alpha\mapsto F^\sigma:=\frac{F}{\prod_{j \notin \sigma} x_j^{\langle m_\sigma,\eta_j\rangle+a_j}}\in\C[X_\sigma].
\end{equation}
Note that with this dehomogenization, the term of $F$ corresponding to $m_\sigma$ becomes constant.
If we take the dual cone $\sigma^{\vee}\subset M\otimes\R$,  we have a natural identification (cf. \cite[Lemma 2.2]{cox95})
\begin{equation}\label{csigma}
\C[X_\sigma]\cong\C[\sigma^\vee\cap M].
\end{equation}

Recall that an $n$-dimensional complex variety $X$ is a $V$-manifold if and only if for every $\bfxi\in X$ there exists a {\em small} finite subgroup $G\subset GL(n,\C)$ such that for some neighborhood 
$W$ of $\bfxi\in X$, we have  that  $(W, \bfxi)\cong(U/G, {\bf0})$, where $U$ is a $G$-invariant neighborhood of ${\bf0}\in\C^n.$ The group $G$ is small if no $g\in G$ has $1$ as an eigenvalue of multiplicity $n-1$, and is unique up to conjugacy.

A simplicial toric variety $X$ is an example of a V-manifold. Indeed, the affine open sets $X_\sigma$ which cover $X$  satisfy 
\begin{equation}\label{coc}
X_\sigma\simeq\C^n/G(\sigma),
\end{equation} 
 with $G(\sigma)$ being  defined as follows: assume w.l.o.g. that $\sigma$ is generated by $\eta_1,\ldots, \eta_n.$ Then, we consider the following exact sequence of $\Z$-modules:
\begin{equation}\label{xxx}
0\to M\stackrel{\gamma}{\to}\Z^n\to D(\sigma)\to0,
\end{equation}
with $\gamma(m):=(\langle m,\eta_1\rangle,\ldots, \langle m,\eta_n\rangle).$ Note that $D(\sigma)$ is finite as $\sigma$ is simplicial. We then set 
\begin{equation}\label{gsigma}
G(\sigma):=\Hom(D(\sigma),\C^*),
\end{equation}
which  is also a finite group.  The map $\Z^n\mapsto D(\sigma)$ induces an action of $G(\sigma)$ to $\C^n,$ and it follows from \cite{cox95} that $X_\sigma$ is the quotient of this action. One can also verify (see for instance the proof of Proposition $3.5$ in \cite{BC94}) that $G(\sigma)$ is small.

If the cone $\sigma$ is simplicial, one can replace the representation of $\C[X_\sigma]$ given in  \eqref{csigma} with an equivalent suitable one. Indeed, thanks to \eqref{coc}, $X_\sigma$ is a quotient of $\C^n$ via the action of the finite commutative group $G(\sigma)$ defined in \eqref{gsigma}, as explained in \cite[Example 1.3.20]{CLS11}. To be more precise, Assume w.l.o.g. that the generators of $\sigma$ are $\eta_1,\ldots, \eta_n,$  and denote by $N_\sigma\subset\Z^n$ the lattice they generate, and by $M_\sigma$ its dual lattice.  As $N_\sigma\subset N,$ we  have $M_\sigma\supset M$, and (Proposition $1.3.18$ in \cite{CLS11})
\begin{equation}\label{rev}
 \C[\sigma^\vee\cap M_\sigma]^{G(\sigma)}=\C[\sigma^\vee\cap M]=\C[X_\sigma],
\end{equation}
In addition, as $\sigma$ is smooth relative to $N_\sigma$ (cf. Example $1.3.20$ in \cite{CLS11}), we deduce that $\C[\sigma^\vee\cap M_\sigma]$ is the coordinate ring of $\C^n$, i.e. it is isomorphic to $\C[u_1,\ldots, u_n],$ a polynomial ring in $n$ variables over $\C.$ So, we then have that 
\begin{equation}\label{zimp}
\C[X_\sigma]\cong \C[u_1,\ldots, u_n]^{G(\sigma)}.
\end{equation}
This identification will be used in the text.

A Cartier  divisor $D$  has no base points, i.e. $\cO(D)$ is generated by its global sections if and only if for all $n$-dimensional $\sigma\in\Sigma,\, m_\sigma\in P_\alpha$
(\cite[Proposition 6.11]{CLS11}).
Sometimes these divisors with no base points are also called semi-ample (cf. \cite{KS05}).
A {\em very ample} divisor $D$ (\cite[Definition 6.1.9]{CLS11}) is one having no base points and also such that the map from $X$ to  $\P(H^0(X,\cO_X(D))^\vee)$ sending a point $\bfxi \in X$ to $(\bfchi^m(\bfxi))_{m\in P_D\cap M}$,
is a closed embedding.
An {\em ample} divisor $D$ is one such that $k D$ is very ample for some $k>0.$   An ample divisor $D$ is very ample if it moreover  satisfies an arithmetic condition regarding the lattice points in the cone bounding $P_D$ from each vertex 

We next relate these notions to the corresponding polytopes.

\subsection{Divisors and polytopes} \label{ssec:divpol}
We highlight the relation between properties of Cartier divisors and properties of the associated polytope. We keep the previous notations. Given a divisor $D$, we denote its class in $A_{n-1}(X)$ by $[D]= \alpha$. All the references in this subsection refer to the book~\cite{CLS11}.

\begin{proposition} \label{prop:nef}
 A divisor $D$ is Cartier and nef if and only if the polytope $P_D$ in \eqref{pd} is a lattice polytope and its vertices are the lattice
  points $m_\sigma$ with $\sigma \in \Sigma(n)$.
\end{proposition}

 The notions of nef and basepoint free divisors are equivalent by Theorem 6.3.12.   Since the fan $\Sigma$ is complete, Proposition~\ref{prop:nef} follows from Theorem 6.1.7.

\begin{proposition} \label{prop:ample}
A Cartier divisor $D$ is ample if and only if it is nef and moreover, $m_\sigma \neq m_{\sigma'}$ for any two different maximal 
dimensional cones $\sigma, \sigma \in \Sigma(n)$. 
\end{proposition}

Proposition~\ref{prop:ample} his follows from Theorem 6.1.7, Lemma 6.1.13 and Theorem 6.1.14.  

\smallskip

We next state the full translation between the fan $\Sigma$ and and the polytope associated with an ample divisor.

\begin{proposition}\label{prop:amplefull}
A nef Cartier divisor $D$  is ample if and only if  $\Sigma$ is the normal fan of $P_D$. 
More explicitly, a nef Cartier divisor $D$ is ample if and only if $P_D$ is a full dimensional  lattice polytope having $n+r$  
facets orthogonal to each of the rays in $\Sigma$, and as many vertices as maximal cones in $\Sigma(n)$,  such that 
the following happens: let $m_\sigma$ be the vertex of $P_D$ corresponding to the cone $\sigma \in \Sigma(n)$ 
and assume $\eta_{i1},\dots, \eta_{ik}$ are the $1$-dimensional cones in $\sigma$; then, there are exactly 
$k$ facets  concurring at $m_\sigma$ which are orthogonal to $\eta_{i1}, \dots, \eta_{ik}$ 
respectively. In particular, $P$ is simplicial iff $\Sigma$ is simplicial. 
\end{proposition}

Proposition~\ref{prop:amplefull} follows from Theorem 6.2.1.

\smallskip

We will also need the following result.

\begin{proposition}\label{prop:1}
If $X$ is a complete toric variety, then the torus $T$ is equal to the intersection of the open sets $X_\sigma$, for all maximal dimensional cones $\sigma$ in the fan of $X$.
\end{proposition}

\begin{proof} As the torus $T$ is contained in any $X_\sigma$, it is contained in their intersection.  The converse
 is a consequence of the Orbit-Cone Correspondence (Theorem 3.2.6). A cone $\sigma$ in $\Sigma(n)$ gives an orbit $O(\sigma)$ and an affine open subset $ X_\sigma$, and by this theorem,
$X_\sigma = \cup_{\tau \text{ face of } \sigma} O(\tau).$ Since orbits are equal or disjoint, it follows that 
\[\cap_{\sigma \in \Sigma} X_\sigma =  \cup_{\tau  \text{ face of every } \sigma \in \Sigma} O(\tau).\]
Now, the fact that $\Sigma$ is complete implies that the zero cone is the only cone that is a face of every cone in $\Sigma$.
\end{proof}

The statement of Proposition~\ref{prop:1} is not true for instance for $X = \mathbb C^n$. There is a single cone $\sigma = \langle e_1, \dots, e_n \rangle$  of maximal dimension and thus the intersection over all the maximal cones equals 
$X_\sigma = \mathbb{C}^n$  and not the torus. Also, the converse of Proposition~\ref{prop:1} does not hold.
 For an easy counterexample in the plane, consider the fan given by the first and third quadrants, denoted by $\sigma_1$ and $\sigma_2$ plus their faces (the four half coordinate axes and the origin).  The associated toric variety is not complete but $X_{\sigma_1} \cap X_{\sigma_2} = T$.

\subsection{Toric varieties from lattice polytopes}\label{ssec:vtp}
Assume that we have lattice polytopes $P_1,\ldots, P_n$ and $P = P_1 + \dots + P_n$ of dimension $n$ as in the Introduction. We can associate an $n$-dimensional complete toric variety to $P$ as follows.
Let $\eta_1,\ldots, \eta_{n+r}$ denote the primitive vectors in $\Z^n$ orthogonal to the $n+r$ facets of $P$. They are the rays of the {\em normal fan} $\Sigma(P)$ of $P$. We  denote by $X_P$ the associated 
complete toric variety $X_{\Sigma(P)}$.

There are integer numbers $a_1, \dots, a_{n+r}$ such that
\begin{equation}\label{P}
P=\{\bfm\in\R^n:\, \langle \bfm,\eta_j\rangle\geq-a_j,\, 1\leq j\leq n+r\}\subset\R^n.
\end{equation}
In the notation of~\eqref{pd}, $P= P_D$, for
$[D] = \sum_{j=1}^{n+r} a_j [D_j]$, with $D_j$ the divisor associated with $\eta_j$.
 There are also integers $a_{ij}$, $1\leq i\leq n$,  $1\leq j\leq n+r$, such that $a_j=a_{1j}+\ldots+a_{nj},$ and 
\begin{equation}\label{PPi}
P_i =\{\bfm \in \R^n\, : \, \langle \bfm,\eta_j \rangle \ge - a_{ij}\,\, 1\leq j\leq n+r\},\ i=1,\ldots, n.
\end{equation}
  As before we denote by $S = \C[x_1,\ldots, x_{n+r}] $ the homogeneous coordinate ring of $X_P$.
  In this ring, we have that $\alpha_{P_i}:=[D_{P_i}]=\sum_{j=1}^{n+r}a_{ij}[D_j]$ for $i=1,\ldots, n$. We also denote  $\alpha_j:=[D_j]$ for $j=1, \dots, n+r$, 
$[D_P]=\sum_{i=1}^n [D_{P_i}]$ and $\alpha_P:=\sum_{i=1}^n \alpha_{P_i}$.  We have that $D_P$ is ample in $X_P$ although not necessarily all the $D_{P_i}$'s must be. They are semi-ample. Each $f_i$ can be then homogeneized to $F_i\in S_{\alpha_{P_i}}$ via \eqref{eFe}

\section{The case of $n+1$ ample divisors}\label{sec:critical}

Suppose that $D_0,\ldots, D_n$ are ample divisors in a complete toric variety $X$ of respective degrees $\alpha_0,\ldots, \alpha_n.$
We define the {\em critical degree} of the sequence $\bfF = (F_0, \dots, F_n)$ of homogeneous polynomials of degrees $(\alpha_0, \dots,\alpha_n)$  as
\begin{equation}\label{ro}
\rho_\bfF:=\left(\sum_{i=0}^n\alpha_i\right)-\rho_0\in A_{n-1}(X).
\end{equation}
In this section, we recall the definition of special homogeneous polynomials of critical degree and we present Proposition~\ref{cd1} on the codimension of the quotient of the homogeneous coordinate ring $S$ of $X$ by the ideal generated by $F_0, \dots, F_n$ at critical degree. We then show in Proposition~\ref{bound} that
controlling the codimension of some graded pieces of $\langle F_0,\ldots, F_n\rangle$  allows us to get  some useful bounds of the dimension of crucial pieces of $S/\langle F_1,\ldots, F_n\rangle$.

\subsection{Special polynomials of critical degree}
As shown in \cite{CCD97}, for every simplicial  $n$-dimensional cone $\sigma\in \Sigma$, if w.l.o.g. $\eta_1,\ldots, \eta_n$ are its generators, every $F_j\in S_{\alpha_j}$ can be written in the form
\begin{equation}\label{desc}
F_j=\sum_{i=1}^n F_{ij}x_i+F_{0j}x_{n+1}\ldots x_{n+r},
\end{equation}
with $F_{ij}\in S$ for all $i,j=0,\ldots, n.$
We set
\begin{equation}\label{delta}
\Delta_{\bfF,\sigma}:=\det(F_{ij})_{0\leq i,j\leq n}\in S_{\rho_\bfF}.
\end{equation} 
 In \cite{DK05}, a generalization of  $\Delta_{\bfF,\sigma}$ was given if there are no simplicial cones.  In this case, $\bfsigma$ actually denotes a flag of cones of $\Sigma: \ \{0\}=\sigma_0\subset\sigma_1\subset\ldots\subset\sigma_n$, with $\dim(\sigma_i)=i.$ For $i=1,\ldots, n,$ we set $z_i$ as the product of all $x_j$ with $\eta_j\in \sigma_j\setminus\sigma_{i-1},$ and then $z_{0}$ being equal to the product of all $x_j$ with $\eta_j\notin\sigma_n.$ Then, as in \eqref{desc}, we have a decomposition of the form
 $F_j=\sum_{i=0}^{n} F_{ij}z_i,$ and with this data we set $\Delta_{\bfF,\bfsigma}$ as the determinant of \eqref{delta}. This element also belongs to $S_{\rho_\bfF}.$

 In another direction, if all the $\alpha_i$'s are equal to a fixed ample degree $\alpha,$ then (cf. \cite[Proposition 4.1]{cox95}) there exists $J_\bfF\in S_{\rho_\bfF}$ such that
\begin{equation}\label{eule}
\sum_{i=0}^n(-1)^i F_i\,dF_0\wedge\ldots\wedge\widehat{dF_i}\wedge\ldots\wedge dF_n= J_\bfF\,\, \Omega,
\end{equation}
with $\Omega$ equal to the {\em Euler form} of $X$ defined as
\begin{equation}\label{omega}
\Omega:=\sum_{I\subset\{1,\ldots, n+r\},\,|I|=n} d_I\prod_{j\notin I} x_j\,d\bfx_I,
\end{equation}
where, if $I=\{i_1,\ldots, i_n\},\, d\bfx_I$ stands for $dx_{i_1}\wedge\ldots\wedge dx_{j_n},$ and $d_I:=\det(\langle m_i,\eta_{i_j}\rangle)_{1\leq i,j\leq n}.$  Note that the product $d_I\,d\bfx_I$ does not depend on the order of the sequence $\{i_1, \dots, i_n\}$.
The polynomial $J_\bfF\in S_{\rho_\bfF}$ is called the {\em toric jacobian} of the system $F_0,\ldots, F_n$. 

This construction has been generalized in \cite{CDS98} to polynomials  of degrees $k_i\alpha,\, i=0,\ldots, n,\, k_i\in\N$, for an ample degree $\alpha$.
In this case, the toric jacobian can be computed as follows (cf. \cite[(1.7)]{CDS98}): if $I=\{i_1,\ldots, i_n\}$ is such that 
$\eta_{i_1},\ldots, \eta_{i_n}$ are linearly independent, then
\begin{equation}\label{JF}
J_\bfF=\frac{1}{d_I\,\prod_{j\notin I}^n x_{j}}\det\left(\begin{array}{cccc}
k_0F_0& k_1F_1&\ldots & k_nF_n\\
\frac{\partial F_0}{\partial x_{i_1}} & \frac{\partial F_1}{\partial x_{i_1}}& \ldots &\frac{\partial F_n}{\partial x_{i_1}} \\
\vdots &\vdots& \vdots&\vdots\\
\frac{\partial F_0}{\partial x_{i_n}}&\frac{\partial F_1}{\partial x_{i_n}}& \ldots & \frac{\partial F_n}{\partial x_{i_n}}
\end{array}\right).
\end{equation}
In general, both the toric jacobian $J_\bfF$ and $\Delta_{\bfF,\bfsigma}$  are  not in the ideal $\langle F_0,\ldots, F_n\rangle,$ but ``almost'' there as the following result shows. 

\begin{proposition}\label{oxi}
For any maximal simplicial cone $\sigma \in\Sigma,$ both $J_\bfF$ and $\Delta_{\bfF,\bfsigma}$ belong to the ideal generated by $F_0^\sigma,\ldots F_ n^\sigma$ in $S_{\bfx^{\widehat{\sigma}}}.$
\end{proposition}

\begin{proof}
Denote by $\{\eta_{i_1},\ldots, \eta_{i_n}\}$ the generators of $\sigma$, and set $I=\{i_1,\ldots, i_n\}.$ Then, $\bfx^{\widehat{\sigma}}=\prod_{j\notin I}^n x_{j}$ and the claim follows straightforwardly for $J_\bfF$ by developing the determinant in \eqref{JF} through the first row of that matrix.

For $\Delta_{\bfF,\bfsigma},$ if there is $i\in\{1,\ldots, n\}$ such that the cone generated by the elements in $\sigma_i\setminus\sigma_{i-1}$ only intersects $\sigma$ in the origin, then $z_i$ is invertible in  $S_{\bfx^{\widehat{\sigma}}},$ and the claim follows straightforwardly by multiplying by $z_i$ the $i$-th row of the matrix \eqref{delta}, and then adding to this row some multiples of the other rows to get $F_0,\ldots, F_n$ appear there, using that $F_\ell=\sum_{j=0}^{n} F_{j\ell}z_j.$ Then, we develop  the determinant via this row and the claim holds. If there is no $i\in\{1,\ldots, n\}$  as above, then each of the $\eta_{i_1},\ldots, \eta_{i_n}$ must belong to one and only one of $\sigma_i\setminus\sigma_{i-1},\, i=1,\ldots, n.$ This implies straightforwardly that  $\sigma=\sigma_n$ and hence $z_0=\bfx^{\widehat{\sigma}}$ is then clearly invertible in $S_{\bfx^{\widehat{\sigma}}},$ and the claim now follows by multiplying by $z_0$ the $0$-th row of \eqref{delta} and proceeding then as before. This concludes with the proof.
\end{proof}

The interest of Proposition \ref{oxi} lies in the fact that this membership is independent of the fact that $V_X(F_0,\ldots, F_n)$ is empty or not. Otherwise the result is trival:
\begin{proposition}\label{goxi}
Let $F_0,\ldots, F_n\in S$ homogeneous such that $V_X(F_0,\ldots, F_n)=\emptyset$.  Then for any maximal cone $\sigma\in\Sigma$, the generated ideal $\langle F_0^\sigma,\ldots F_ n^\tau \rangle$ in $S_{\bfx^{\widehat{\sigma}}}$ is equal to the whole ring $S_{\bfx^{\widehat{\sigma}}}$.
\end{proposition}

\begin{proof}
From \cite[Prop.~5.2.8]{CLS11}, we deduce that  the affine variety $V_{\C^{n+r}}(F_0,\ldots, F_n)$ is contained in the variety defined by the irrelevant ideal $\langle \bfx^{\widehat{\sigma}}, \sigma\in\Sigma\rangle.$ Then,  there must exist $m\in\N$ such that  $(\bfx^{\widehat{\sigma}})^m\in\langle F_0,\ldots, F_n\rangle.$ Hence, $(\bfx^{\widehat{\sigma}})^mG$ belongs to this ideal for any homogeneous polynomial $G$ and the claim follows straightforwardly.
\end{proof} 

The following result can be found in \cite[Theorems 0.2 and 0.3]{CCD97}, \cite[Theorem 5.1]{cox96}, and \cite[Theorem 2.1]{CD05}.
\begin{theorem}\label{if} Let $X$ be a complete  toric variety, and assume that $F_i\in S_{\alpha_i}$ for $i=0,\ldots, n,$ where each $\alpha_i$ is  ample, and  the $F_i$'s do not vanish simultaneously on $X$. \begin{enumerate}
\item If $\bfsigma$ is a flag of cones in $\Sigma,$ then 
 $\Delta_{\bfF,\bfsigma}\notin \langle F_0,\ldots, F_n\rangle.$
\item If all the $\alpha_i=k_i\alpha$ with $\alpha$ ample, then  $J_\bfF\notin \langle F_0,\ldots, F_n\rangle.$
\end{enumerate}
\end{theorem}

In general, not all the divisors associated to polytopes are ample, but semi-ample. In \cite{KS05}, some determinants like \eqref{delta} are considered associated to colouring partitions of the input supports.

\subsection{Codimension one}\label{ssec:cod1}

To  put Theorem \ref{if} within context, in \cite{DK05}, when all the divisors $\alpha_i$ are ample the following map was considered (3.7) y en (3.8):
\begin{equation}\label{phiD}
\begin{array}{ccc}
\phi_{\bfF,\sigma}:S_{\rho_\bfF-\alpha_0}\oplus\ldots\oplus S_{\rho_\bfF-\alpha_n}\oplus S_0&\to& S_{\rho_\bfF}\\
(G_0,\ldots, G_n,\lambda)&\mapsto & \sum_{i=0}^n G_i F_i +\lambda\Delta_{\bfF,\bfsigma},
\end{array}
\end{equation}
Previously, in \cite{CDS98}, the following variation of \eqref{phiD} was studied:
\begin{equation}\label{phiJ}
\begin{array}{ccc}
\phi_{J_\bfF}:S_{\rho_\bfF-k_0\alpha}\oplus\ldots\oplus S_{\rho_\bfF-k_n\alpha}\oplus S_0&\to& S_{\rho_\bfF}\\
(G_0,\ldots, G_n,\lambda)&\mapsto & \sum_{i=0}^n G_i F_i +\lambda J_\bfF,
\end{array}
\end{equation}
with $k_0,\ldots, k_n\in\N,\,\alpha$ being an ample divisor, and $\alpha_i=k_0\alpha,\,i=0,\ldots, n.$

We then have the following result.
\begin{theorem}\label{resn0} (\cite[Prop 2.1]{CDS98}, \cite[Prop 4.2]{DK05})
Let $F_0,\ldots, F_n\in S$ be homogeneous of ample degress $\alpha_0,\ldots, \alpha_n$ as in \eqref{phiD} (resp. \eqref{phiJ}). The map $\phi_{\bfF,\sigma}$ (resp. $\phi_{J_\bfF}$) is onto if and only if $V_X(F_0,\ldots, F_n)=\emptyset.$ In particular, $\dim_\C\left(S/\langle F_0,\ldots, F_n\rangle\right)_{\rho_\bfF}=1$ if this is the case.
\end{theorem}

Theorem~\ref{fonlyif} implies that if any of these maps fail to be onto is because either $\Delta_{\bfF,\bfsigma}$ in the case of $\phi_{\bfF,\sigma}$ or the jacobian $J_\bfF$ in the caso of $\phi_{J_\bfF}$ belong to the ideal $\langle F_0,\ldots, F_n\rangle$.  From a computational perspective, testing this latter property is more efficient than checking if the whole linear map \eqref{phiD} or \eqref{phiJ} is onto. 
As mentioned above, the validity  of Theorem~\ref{fonlyif} is known for the case of $X=\P^n$. Indeed, an elementary proof of this result has been given in \cite{spo89}, and extensions to affine and local scenarios have been explored in \cite{BY00} and \cite{hic08} respectively.  

\smallskip

One can generalize the codimension one property given in Theorem \ref{resn0} above to more general degrees than being all of them ample.  We extend Definition~\ref{def:polytopes} to a sequence of $n+1$ lattice polytopes with the same conditions for any subset of up to $n$ of them for the property of being an {\em essential } sequence and up to $n-1$ of them for the property of being an {\em indecomposable} sequence. 

 Theorem $2.2$ in \cite{CD05} shows that the codimension one property does not hold for any sequence of nef Cartier divisors, but it implies the following result for indecomposable sequences.

\begin{proposition}\label{cd1}
Assume that $\alpha_0,\ldots, \alpha_n$ are degrees of nef Cartier divisors  such that the sequence of corresponding polytopes $P_{\alpha_0}, \dots, P_{\alpha_n}$ is indecomposable. 
Consider $F_0,\ldots, F_n\in S$, homogeneous  of respective degrees $\alpha_0,\ldots, \alpha_n$ such that $V_X(F_0,\ldots, F_n)=\emptyset.$ Then, $\dim_\C\left(S/\langle F_0,\ldots, F_n\rangle\right)_{\rho_\bfF}=1.$
\end{proposition}

We end this section with a results that relates the dimension of quotients by $n$ and $n+1$ divisors.

\begin{proposition}\label{bound}
Assume that  $\alpha_0,\ldots, \alpha_n$ are degrees of nef Cartier divisors, with $\alpha_0$ being also ample, satisfying  $\alpha_0\geq \rho.$ If $F_0,\ldots, F_n\in S$ are homogeneous  of respective degrees $\alpha_0,\ldots, \alpha_n$ such that $V_X(F_0,\ldots, F_n)=\emptyset,$ then  
$$\dim_\C\left(S/\langle F_1,\ldots, F_n\rangle\right)_{\alpha_1+\ldots+\alpha_n-\rho_0}\geq \deg\big(V_X(F_1,\ldots, F_n)\big)- \dim_\C\left(S/\langle F_0,\ldots, F_n\rangle\right)_{\rho_\bfF}.
$$
\end{proposition}
\begin{proof}
As $X$ is projective (because $\alpha_0$ is ample), the fact that $V_X(F_0,\ldots, F_n)=\emptyset$ implies  $V_X(F_1,\ldots, F_n)$ being zero-dimensional follows from \cite[Lemma $4.14$]{CCD97}.
Consider the following exact sequence of $\C$-vector spaces:

$$0\to K\to  S/\langle F_1,\ldots, F_n\rangle_{\alpha_1+\ldots+\alpha_n-\rho_0}\stackrel{\phi_1}{\to}S/\langle F_1,\ldots, F_n\rangle_{\rho_\bfF}\stackrel{\phi_2}{\to} S/\langle F_0,\ldots, F_n\rangle_{\rho_\bfF}\to0,$$
with $\phi_1$ being the multiplication by $F_0,\, \phi_2$ the projection to the quotient, and $K=\mbox{Ker}(\phi_1).$  From \cite[Theorem $4.4$]{BT22}, due to the hypothesis imposed on $\alpha_0,$ we have that $\dim_\C\left(S/\langle F_1,\ldots, F_n\rangle_{\rho_\bfF}\right)= \deg\big(V_X(F_1,\ldots, F_n)\big).$ By counting dimensions of all these spaces in the above sequence, the claim then follows straightforwardly. 
\end{proof}
 
\bigskip
\section{Toric residues and residues in the torus}\label{sec:3}
Toric jacobians  are closely related with the computation of toric residues, which were introduced and studied in \cite{cox96, CCD97}.  We recall their definition and properties in this section. For simplicity and historical reasons, we present all of our results for the case where the ground field is $\C,$ but these results are actually valid for any algebraically closed field of characteristic zero, see for instance Chapter $14$ of \cite{kun08}. 
We first recall the notion of global residues in section~\ref{ssec:gr}. In section~\ref{lrd} we summarize local residues and local duality in the smooth case. Section~\ref{groth} recalls the definition of local residues in $V$-manifolds and Section~\ref{tgr} deals with toric and global residues in the simplicial case. 

\subsection{Global residues}\label{ssec:gr}
For a sequence of polynomials $F_0,\ldots, F_n\in S$ of degrees $\alpha_0,\ldots, \alpha_n$, every $H\in S_{\rho_\bfF}$ (where $\rho_\bfF$ has been defined in \eqref{ro}) determines a meromorphic $n$-form in $X$
\begin{equation}\label{omegaf}
\omega_\bfF(H):=\frac{H\Omega}{F_0\ldots F_n}
\end{equation}
with $\Omega$ defined in \eqref{omega}. If the $F_i$'s do not vanish simultaneously in $X,$ then relative to the open cover $U_i:=\{\bfxi \in X:\, F_i(\bfxi)\neq0\}$ of $X,$ \eqref{omegaf} gives a C$\breve{e}$ch cohomology class $[\omega_\bfF(H)]\in H^n(X,\widehat{\Omega}^n_X),$ with $\widehat{\Omega}^n_X$ being the sheaf of Zariski $n$-forms on $X,$ (i.e. $\widehat{\Omega}^n_X=j_*\Omega^n_{X_0},$ where $X_0$ is the smooth part of $X$ and $j:X_0\hookrightarrow X$ is the natural inclusion). One can check that  if $H\in\langle F_0,\ldots, F_n\rangle,$ then $\omega_\bfF(H)$ is a C$\breve{e}$ch coboundary. Thus, $[\omega_\bfF(H)]$ depends only on the equivalence class of $H$ modulo $\langle F_0,\ldots, F_n\rangle.$ The {\em toric residue} map is then defined as follows:
\begin{equation}\label{tres}
\begin{array}{cccc}
\Res_\bfF:&S_{\rho_\bfF}/\langle F_0,\ldots, F_n\rangle_{\rho_\bfF}&\to&\C\\
&[H]&\mapsto & \mbox{Tr}_X([\omega_\bfF(H)]),
\end{array}
\end{equation}
where $\mbox{Tr}_X:H^n(X, \widehat{\Omega}^n_X)\to\C$ is the trace map. We will simply denote it by $\Res_\bfF(H)$ to abuse notation.

The following elements have non trivial residue.
\begin{proposition}\label{deta}(\cite[Theorem 0.2]{CCD97} and \cite[Theorem 2.1]{DK05})
If $X$ is complete, $\bfsigma$ a flag of $n+1$ cones in $\Sigma,$ each $\alpha_i$ ample and the $F_i$'s do not vanish simultaneously in $X,$ then
$$\Res_\bfF(\Delta_{\bfF,\bfsigma})=\pm1.
$$
\end{proposition}
For a divisor $D$ of degree $\alpha,$ Recall the definition of the polytope $P_D$ given in \eqref{pd}. 
Its  normalized volume in $M\otimes\R$ (where the fundamental simplex in $M$ has volume $1$) will be denoted with $\mbox{vol}(P_D).$
\begin{proposition}\label{jaobian}(\cite[Theorem 1]{CD97})
If $X$ is a complete toric variety, $\alpha\in A_{n-1}(X)$ an ample degree, and each $\alpha_i=k_i\alpha$ for positive integers $k_0,\ldots, k_n$, then
\begin{equation}\label{D}
\Res_\bfF(J_\bfF)=\Big(\prod_{i=0}^nk_i\Big)\vol(P_D).
\end{equation}
\end{proposition}

\begin{remark}\label{desing}
If $\Sigma'$ is a refinement of $\Sigma$ which does not introduce new rays, then it turns out that both toric varieties $X_{\Sigma'}$ and $X_\Sigma$ share the same coordinate ring and, as it was shown in  \cite[Proposition $2.1$]{KS05},
$F_0,\ldots, F_n$ does not have any common zeroes in $X_{\Sigma'}$ if this happens in $X_\Sigma.$ Moreover,  $\mbox{Tr}_{X_\Sigma}([\omega_\bfF(H)])=  \mbox{Tr}_{X_\Sigma'}([\omega_\bfF(H)])$ for any homogeneous $H\in S_{\rho_\bfF}.$ 
\end{remark}

\medskip
\subsection{Local residues and duality in the smooth case}\label{lrd}
As explained in \cite{DS91} and in the introduction of \cite{CDS94}, the local residue is a generalization of the Cauchy formula in several complex variables. We recall here its analytic definition and main properties.
For a point $\bfxi\in\C^n,$ denote with $\cO_\bfxi$ the stalk of analytic functions over $\bfxi.$ Let $f_{1\bfxi},\ldots, f_{n\bfxi}\in\cO_\bfxi,$ and $U$ a ball centered at $\bfxi$ such that there exist representatives $f_1,\ldots, f_n\in\cO(\overline{U})$ verifying $V_{\overline{U}}(f_1,\ldots,f_n)=\{\bfxi\}.$ For any meromorphic $n$-form $\omega$ in $U$ having its poles on $\cup_{i=1}^n (f_i=0),$ the residue of $\omega$ at $\bfxi$ is the complex number 
\begin{equation}\label{reslocal}
\mbox{Res}_\bfxi(\omega):=\frac{1}{(2\pi i)^n}\int_{\{\bfz\in U:\,|f_j(z)|=\delta_j,\,1\leq j\leq n\}}\omega,
\end{equation}
with $\delta_j>0$ being such that the set $\{(w_1,\ldots, w_n)\in\C^n:\,|w_j|=\delta_j,\, 1\leq j\leq n\}$ is contained in the open set $\bff(U),$ and $\bff:=(f_1,\ldots, f_n),$ and satisfying $$\{\bfz\in\overline{U}:\,|f_j(\bfz)|=\delta_j,\, 1\leq j\leq n\}\cap\partial U=\emptyset.$$ The tube $\{\bfz\in u:\,|f_j(\bfz)|=\delta_j,\,1\leq j\leq n\}$ is compact, of real dimension $n,$ and  the integral is computed with the orientation given by the form $d(\arg f_1)\wedge\ldots\wedge d(\arg f_n).$

It can be shown that this definition is independent of the neighborhood $U$ chosen, the values of $\delta_j,\,1\leq j\leq n,$ and also of the representatives $f_1,\ldots, f_n\in\cO(\overline{U})$ of $f_{1\bfxi},\ldots, f_{n\bfxi}.$

The local residue allows us to define a $\C$-linear operator 
\begin{equation}\label{resop}
\begin{array}{cccc}
R_{\bff, \bfxi}[\omega]:& \cO_\bfxi& \to& \C\\
  &  g& \mapsto & \mbox{Res}_\bfxi(g\cdot\omega).
\end{array}
\end{equation}
In the case $\omega=\frac{h\,dz}{f_1\ldots f_n},$  the following property of the local residue will be of use for us in the sequel:
\begin{proposition}\label{lresid}
Let $h$ be such that $h_\bfxi\in\cO_\bfxi,$ Then,
 $h_\bfxi\in I_\bfxi\iff R_{\bff, \bfxi}\big[\frac{h\,dz}{f_1\ldots f_n}\big]\equiv0,$ the zero operator. 
\end{proposition}

If we pick $h= J_\bff,$ the jacobian of the system $f_1,\ldots, f_n,$ it can be shown that $\mbox{Res}_\bfxi(\omega)=\mu_{\bff,\bfxi},$ the local multiplicity of $\bfxi$ as a zero of $\bff.$

\medskip
\subsection{Residues in simplicial affine varieties }\label{groth}
We recall here the definition and main properties of \textit{local Grothendieck residues} at a point $\bfxi$ in a \textit{V-manifold} or \textit{orbifold} $X$ extracted from \cite[\S 4]{CCD97}. If $X$ is a simplicial toric variety, this is the case.  Note that as these are local properties and the torus $(\C^*)^n$ is a subvariety of $X$ of the same dimension, we are in particular defining local Grothendieck residues at points $\bfxi\in (\C^*)^n.$

Let $\bfxi$ a point in a V-manifold $X$, and   $(W_1, \bfxi)\cong (U_1/G,{\bf 0})$  be a standard model. Let $W$ be a relatively compact neighborhood of $\bfxi$ such that 
$W\subset\overline{W}\subset W_1,$ and $U$ be a $G$-invariant neighborhood of ${\bf0}$ such that $W\cong U/G$ and $U\subset\overline{U}\subset U_1.$ Suppose that
$f_1,\ldots ,f_n \in\cO(\overline{W})$ have $\bfxi$ as their only common zero in $\overline{W}.$ Pulling-back to $\overline{U},$ it follows 
that the hypersurfaces $\tilde{D}_i = \{\tilde{f}_i = 0\}$ intersect only at ${\bf0}.$ Given now a meromorphic $n$-form $\omega$ on $W$ with polar set contained in $\cup_{i=1}^n \{f_i = 0\},$ we denote by $\tilde{\omega}$ its pull-back to $U$ and define the local Grothendieck residue 
\begin{equation}\label{rloc}
\mbox{Res}_\bfxi(\omega):=\frac{1}{|G|}\mbox{Res}_{\bf0}(\tilde{\omega}),
\end{equation}
where
$\mbox{Res}_{\bf0}(\tilde{\omega})$ is defined as in \eqref{reslocal}  with respect to $\tilde{f}_1,\ldots, \tilde{f}_n.$

The connection between local Grothendieck residues and the toric residue defined above is the following:
\begin{theorem}\cite[Theorem 0.4]{CCD97}\label{gloc}
If $X$ is a complete simplicial toric variety of dimension $n,$ and $F_0,\ldots, F_n\in S$ homogeneous which do not vanish simultaneously on $X,\, H\in S_{\rho_\bfF},$ and 
$V_X(F_1,\ldots, F_n)$ finite, then the $n$-form $\omega_H:=\frac{H/F_0\,\Omega}{F_1\ldots F_n }$ is meromorphic  in $X$ with polar set contained in $\cup_{i=1}^n \{F_i = 0\},$ and moreover
$$\mbox{Res}_\bfF(H)=\sum_{\bfxi\in V_X(F_1,\ldots, F_n)}\mbox{Res}_\bfxi(\omega_H).
$$
\end{theorem}

This result has the following interesting property which will be well used in the sequel:
 \begin{corollary}\label{ejsi}
If $X$ is a complete simplicial toric variety of dimension $n,\, F_1,\ldots, F_n\in S$ homogeneous of degrees $\alpha_1,\ldots, \alpha_n$ respectively be such that  $V_X(F_1,\ldots, F_n)$ is finite, and $G\in S_{\alpha_1+\ldots+\alpha_n-\rho_0}.$ Then 
$$\sum_{\bfxi\in V_X(F_1,\ldots, F_n)}\mbox{Res}_\bfxi\Big(\frac{G\,\Omega}{F_1\ldots F_n}\Big)=0.
$$
\end{corollary}
\begin{proof}
Let $F_0\in S$ be a homogeneous form such that $V_X(F_0,\ldots, F_n)=\emptyset,$ and set $H:=F_0\cdot G.$ The result then follows straightforwardly from Theorem \ref{gloc}.
\end{proof}

\medskip

\subsection{Toric and global residues in the simplicial case}\label{tgr}
Assume that we have lattice polytopes $P_1,\ldots, P_n$ and $P$ as in the Introduction and consider the variety $X_P$ defined in Section~\ref{ssec:vtp}.
Denote by $F_i\in S_{\alpha_{P_i}}$ the $P_i$-th homogenization of $f_i,$ as in \eqref{eFe}. 

Given a homogeneous polynomial $H\in\C[x_1,\ldots,x_{n+r}]$ of critical degree $\rho_\bfF,$ the expression
\begin{equation}\label{uu}
\frac{H\Omega}{F_0\ldots F_n}
\end{equation}
defines a meromorphic $n$-form on $X_P,$ whose restriction to the torus $(\C^*)^n$ may be written via \eqref{tx} as
\begin{equation}\label{dd}
\frac{h}{f_0\ldots f_n}\frac{dt_1}{t_1}\wedge\ldots\wedge\frac{dt_n}{t_n},
\end{equation}
with $h\in\C[\bft^{\pm1}]$ supported in $(P_0+\ldots+P_n)^\circ.$

If $F_0,\ldots, F_n$ do not have common zeroes in $X_P,$  then the function  $h/f_i$ is regular at these points. If $X$ is simplicial, from \eqref{uu} and \eqref{dd} we then have that, for  each $\bfxi\in V_X(F_1,\ldots, F_n)\cap (\C^*)^n,$ 
\begin{equation}\label{sosi}
\mbox{Res}_\bfxi\Big(\frac{G\,\Omega}{F_1\ldots F_n}\Big)=\Res_\bfxi\left(\frac{h}{f_0\ldots f_n}\frac{dt_1}{t_1}\wedge\ldots\wedge\frac{dt_n}{t_n}\right),
\end{equation}
where the left hand side is the residue in the orbifold $X,$ and the right hand side is the local residue in $\C^n.$
In addition, if for some $i\in\{0,\ldots, n\},$ the set $\{x\in X:\,F_j(x)=0\, \forall j\neq i\}$ lies in the torus $(\C^*)^n,$
then the following result holds independently of whether $X$ is simplicial or not:
\begin{proposition}\cite[Proposition 1.3]{CDS98}\label{locglob}
$$\Res_\bfF(H)=(-1)^i \Res^T_{f_0,\ldots, \widehat{f_i},\ldots, f_n}(h/f_i).$$
\end{proposition}

We remark again that all these results can be defined and extended to the case where the base field is any algebraically closed field of characteristic zero, see Chapter $14$ of \cite{kun08}. For instance, Theorem \ref{gloc} above is proven with weaker hypothesis over any field in Theorem $14.18$ in \cite{kun08}.

\bigskip
\section{Residues and local multiplicities of points in affine simplicial toric varieties}\label{sec:4}
We present in this section the approach and results of \cite[Section 3.1]{BT22} on the local rings of points in toric varieties lying in open sets associated with simplicial fans, and we relate them to residues. 

\subsection{Local rings and multiplication maps} \label{ssec:localring}

Let $X$ be a complete toric variety with fan $\Sigma$ as in Section \ref{sec:1}. 
Let $\sigma\in\Sigma$ be an $n$-th dimensional cone, and denote with $m_1,\ldots, m_\ell$ the generators of $\sigma^\vee\cap M$ as an additive semigroup, then it turns out that -thanks to \eqref{csigma}-
\begin{equation}\label{konk}
\C[X_\sigma]\cong\C[y_1,\ldots, y_\ell]/ I_{\cA},
\end{equation} where $I_\cA$ is a toric ideal.

\begin{example}
Consider for instance the cone cone $\sigma$ in dimension $3$, such that the dual cone $\sigma^\vee$ has $1$-dimensional rays generated by the lattice vectors $\eta_1=(1,0,0)$, $\eta_2=(1,2,0)$, and $\eta_3=(0,0,1)$. Then, $\sigma$ is simplicial but it is not smooth because
the determinant of these vectors is equal to $2$. To generate the lattice points in
$\sigma^\vee$, we need to add the vector
$\eta_4=(1,1,0)$. Then, $\ell =4$, and as
$2 \eta_4 = \eta_1 + \eta_2$, we have that~\eqref{konk} takes in this case the explicit form:
\begin{equation}\label{eq:presentation}
\C[X_\sigma]\cong\C[y_1,\ldots, y_4]/ \langle y_4^2- y_1 y_2 \rangle.
\end{equation}
Note that even if $X_\sigma$ is the quotient by a finite group as in~\eqref{coc} (indeed, in this case $G(\sigma)$ is equal to the group of quadratic roots if unity, with $2$ elements $1$ and $-1$), the origin is not
an isolated singularity of $X_\sigma$. Indeed, as the gradient of $y_4^2- y_1 y_2$ vanishes
at any point of the form $(0,0, y_3, 0)$, all these points are singular.
\end{example}

For any ideal $I_\sigma\subset\C[X_\sigma],$ we have a map $\C[y_1,\ldots, y_\ell]\to\C[X_\sigma]\to\C[X_\sigma]/I_\sigma,$ which allows us to consider $V_{X_\sigma}(I_\sigma)$ as either $\mbox{Spec}\big(\C[X_\sigma]/I_\sigma\big)$ or $\mbox{Spec}\big(\C[y_1,\ldots, y_\ell]/I_{\sigma,\cA}\big)$ with $I_\cA\subset I_{\sigma,\cA}\subset\C[y_1,\ldots, y_\ell].$ Note that $\sigma$ is a smooth cone if and only if $\ell=n$ and $I_\cA=\bf0.$

As a consequence, we have that $V_{X_\sigma}(I_\sigma)$ is zero-dimensional if and only if the ideal  $I_{\sigma,\cA}$ defines a zero-dimensional scheme in $\C[y_1,\ldots, y_\ell].$ If this is the case, the multiplicity structure of the points  can be described via the primary decomposition of this ideal and local differential operators as follows: from a minimal primary decomposition
$I_{\sigma,\cA}=Q_1\cap\ldots\cap Q_\delta,$ we get 
\begin{equation}\label{tcr}
\C[y_1,\ldots, y_\ell]/ I_{\sigma,\cA}\cong \C[y_1,\ldots, y_\ell]/Q_1\times\ldots\times\C[y_1,\ldots, y_\ell]/Q_\delta,
\end{equation}
and deduce that $V(I_{\sigma,\cA})$ has $\delta$ points $\bfxi_1,\ldots, \bfxi_\delta\in\C^\ell,$ each of them with multiplicity $\mu_j=\dim_\C\big(\C[y_1,\ldots, y_\ell]/Q_j\big),\,j=1,\ldots, \delta.$
Let $\cD:=\C[\partial_1,\ldots, \partial_\ell]$ be the ring of polynomials in partial derivatives, with $\partial_i:=\frac{\partial}{\partial y_i},\,i=1,\ldots,\ell.$ We define, for $\bfa=(a_1,\ldots, a_\ell)\in\Z_{\geq0}^\ell,$ the element 
$\partial_\bfa:=\frac{1}{a_1!\ldots a_\ell!}\frac{\partial^{a_1+\ldots+a_\ell}}{\partial y_1\ldots \partial y_\ell}\in\cD.$
For $\bfb\in\Z_{\geq0}^\ell,$ we also define the $\C$-linear operator $s_\bfb:\cD\to\cD$ as follows: if 
$$s_\bfb\Big(\sum_\bfa c_\bfa\partial_\bfa\Big):=\sum_{\bfa-\bfb\in\Z_{\geq0}^\ell}c_\bfa\partial_{\bfa-\bfb}.
$$
With the help of this operator, Leibniz'rule for the derivation of a product has an easy expression:
\begin{equation}\label{lbtz}
\partial(g\cdot f)=\sum_{\bfb\in\Z_{\geq0}^\ell} \partial_\bfb(g)\big(s_\bfb(\partial)\big)(f) \ \forall \partial\in\cD,\, \forall g,\,f\in\C[y_1,\ldots, y_\ell].
\end{equation}
A $\C$-vector subspace $V\subset\cD$ is said to be \textit{closed} if $\dim_\C(V)<\infty,$ and for all $\partial\in V,\, \bfb\in\Z_{\geq0}^\ell, s_\bfb(\partial)\in V.$

\begin{proposition}[BT20, page 2407]\label{bfr}
For each of the primary components $Q_i$ from above, $1\leq i\leq\delta,$ there exists a unique closed subspace $V_i\subset\cD$ such that
$$Q_i=\{f\in\C[y_1,\ldots, y_\ell]:\, \partial(f)(\bfxi_i)=0 \ \forall \partial\in V_i\}.
$$
\end{proposition}
For $\bfxi\in\C^\ell$ we denote with $\mbox{ev}_\bfxi:\C[y_1,\ldots, y_\ell]\to\C$ the evaluation morphism $\mbox{ev}_\bfxi(f):=f(\bfxi).$ We then have the following
\begin{corollary}\label{dual}
For $i=1,\ldots, \delta,$
$$\mbox{ev}_{\bfxi_i}\circ V_i:=\{\mbox{ev}_{\bfxi_i}\circ\partial \in(\C[y_1,\ldots,y_\ell]) ^\vee\}]\cong \Big(\C[y_1,\ldots, y_\ell]/ Q_i\Big)^\vee.$$
\end{corollary}
The \textit{order} of a non zero  element $\sum_\bfa c_\bfa\partial_\bfa\in\cD$ is defined as the maximum of $a_1+\ldots+a_\ell$ among those $\bfa=(a_1,\ldots, a_\ell)$ such that $c_\bfa\neq0.$
For a $\C$-vector subspace $V\subset\cD,$ we denote with $V^d$ the subspace of elements of $V$ with order less than or equal to $d.$
A sequence $\partial^1,\ldots, \partial^\mu$ is said to be a \textit{consistently ordered basis} for $V$ if it is a basis of this subspace, and for every $d\geq0,$ there exists $\mu_d\in\N$ such that 
$\partial^1,\ldots, \partial^{\mu_d}$ is a basis of $V^d.$ From this definition it is clear that $\partial^1$ must always be a nonzero constant in such a basis, which we will assume it will be equal to $1.$ 
It is also easy to deduce that any closed subspace has a consistently ordered basis, so we get the following:
\begin{proposition}\label{cun}
For $i=1,\ldots, \delta,$ the subspaces $V_i$ from Proposition \ref{bfr} have a consistenly ordered basis $\{\partial^{i,1},\ldots, \partial^{i,\mu_i}\},\ 1\leq i\leq\delta.$
\end{proposition}
For $i=1,\ldots, \delta$ and $j=1,\ldots, \mu_i$ set  
\begin{equation}\label{pij}
\partial_{i,j}:=\mbox{ev}_{\bfxi_i}\circ\partial^{i,j}\in  (\C[y_1,\ldots, y_\ell]/ Q_i)^\vee.
\end{equation}
From Proposition \ref{cun} and Corollary \ref{dual}, we then have that $\{\partial_{i,1},\ldots, \partial_{i,\mu_i}\}$ is a basis of $ \Big(\C[y_1,\ldots, y_\ell]/ Q_i\Big)^\vee.$ We call this also a \textit{consistently ordered basis} of this subspace.

For $g\in\C[y_1,\ldots, y_\ell],$ consider the multiplication map $M_g:\C[y_1,\ldots, y_\ell]/I_{\sigma,\cA}\to\C[y_1,\ldots, y_\ell]/I_{\sigma,\cA},$ given by the multiplication by the class of $g$ in the quotient. We use the adjoint application of this map $M^{\#}_g:\Big(\C[y_1,\ldots, y_\ell]/I_{\sigma,\cA}\Big)^\vee\to\Big(\C[y_1,\ldots, y_\ell]/I_{\sigma,\cA}\Big)^\vee$ and compute its matrix with respect to the sequence  $\{\partial_{1,1},\ldots, \partial_{\delta,\mu_\delta}\},$ which is a basis thanks to \eqref{tcr}.

To do this, note that for a given $f\in\C[y_1,\ldots, y_\ell],$ by Leibnitz's rule we have that
\begin{equation}\label{dprod}
\partial_{i,j}\big(M_g([f])\big)=\partial_{i,j}(g\cdot f)=\sum_{\bfb\in\N_{\geq0}^\ell} \partial_\bfb(g)s_\bfb(\partial^{i,j})(f) |_{\bfy=\bfxi_i}.
\end{equation}
By using the fact that $\partial^{i,1}$ is equal to the identity, and that if $\bfb\neq\bf0,$ then $s_\bfb(\partial^{i,j})$ belongs to the $\C$-vector space generated by $\partial^{i,1},\ldots, \partial^{i,j-1},$ we then have that there exists complex numbers $c_{i,j,k}$  with $k<j,$ such that
\begin{equation}\label{triang}
\begin{pmatrix}
\partial_{i,1}\\ \partial_{i,2}\\ \vdots\\  \partial_{i,\mu_i}
\end{pmatrix}\circ \, M_g= 
\begin{pmatrix}
g(\bfxi_i) & 0 & \ldots & 0\\
c_{i,2,1}& g(\bfxi_i) &  \ldots &0  \\
\vdots & & \ddots& \vdots \\
c_{i,\mu_i,1}& c_{i,\mu_i,2} & \ldots & g(\bfxi_i)
\end{pmatrix}
\cdot
\begin{pmatrix}
\partial_{i,1}\\ \partial_{i,2}\\ \vdots\\  \partial_{i,\mu_i}
\end{pmatrix}.
\end{equation}

\medskip
\subsection{Local rings and residues in the simplicial case}\label{sss}
If $\sigma$ is a simplicial cone, via the description given in \eqref{rev} and \eqref{zimp}, we can  study the affine primary decomposition of any zero-dimensional ideal contained in $\C[X_\sigma],$ and it is an easy exercise to  verify that for homogeneous polynomials $F_1,\ldots, F_n\in\C[x_1,\ldots, x_{n+r}]$
of degree $[D]$ with $D$ as in \eqref{alpa} with $a_1=\ldots a_n=0$ (i.e. ${\bf0}\in M$ is the vertex of the polytope $P_D$ defined by $\eta_1,\ldots, \eta_n$), and
having a finite number of zeroes in $X_\sigma$, the polynomials $\tilde{F}_1=F_1(u_1,\ldots, u_n,1,\ldots, 1),\ldots, \tilde{F}_n=F_n(u_1,\ldots, u_n, 1,\ldots, 1)$ actually belong to $\C[u_1,\ldots, u_n]^{G(\sigma)}.$ Indeed, if we identify the variables $u_1,\ldots, u_n$ with the dual basis of $\{\eta_1,\ldots, \eta_n\}$ in $N_\sigma$, from \eqref{eFe} and \eqref{rev} we deduce that
$$
\tilde{F}_i=\sum_{m\in P_D\cap M}c_m\prod_{j=1}^{n}u_j^{\langle m,\eta_j\rangle},
$$
and note that $(\langle m,\eta_1\rangle,\ldots,\langle m,\eta_n\rangle)$ is equal to $[{\bf0}]$ in $D(\sigma)$ from \eqref{xxx}, so this polynomial is $G(\sigma)$-invariant.

For any $H\in\C[u_1,\ldots, u_n]^{G(\sigma)}/\langle\tilde{F}_1,\ldots, \tilde{F}_n\rangle,$ 
the meromorphic differential form
\begin{equation}\label{ome}
\frac{H}{\tilde{F}_1\ldots \tilde{F}_n}\,\frac{du_1}{u_1}\wedge\ldots\wedge\frac{du_n}{u_n}
\end{equation}
is $G(\sigma)$-invariant. If in addition ${H}$ is a multiple of $u_1\ldots u_n,$ then its polar set is contained in $\cup_{i=1}^n\{\tilde{F}_i=0\}.$  Assume further that ${\bf0}\in\C^n$ is an isolated zero of $V_{\C^n}(\tilde{F}_1,\ldots, \tilde{F}_n)$. We can then apply to \eqref{ome} the local residue \eqref{rloc} in the orbifold $X_\sigma.$  Denote with $\bfxi\in X_\sigma$ the point whose representative in $\C^n/G(\sigma)$ is ${\bf0},$ and define then the following operator $\mbox{Res}_{F,\bfxi}: \C[u_1,\ldots, u_n]^{G(\sigma)} \to \C$
\begin{equation}\label{eq:resop}
\begin{array}{cccl}
 H & \mapsto & \left\{\begin{array}{cc}
\frac{1}{|G(\sigma)|}\mbox{Res}_{0} \left(\frac{H}{\tilde{F}_1\ldots \tilde{F}_n}\,\frac{du_1}{u_1}\wedge\ldots\wedge\frac{du_n}{u_n}\right) & \mbox{if} \,
u_1\ldots u_n \, | \, H\\
0 & \mbox{if not.}
\end{array}
 \right.
\end{array}
\end{equation}
Note that the nontrivial formula above  is the local residue $\mbox{Res}_\bfxi$ defined in \eqref{rloc} of the diferential form whose pushback to $\C^n$ is \eqref{ome}.

We will now relate  this construction with the dual of $\C[X_\sigma]/I_\sigma $ presented in \eqref{tcr}.
We assume  that
$V_{X_\sigma}(I_\sigma)$ is zero-dimensional, and denote as before $\bfxi_1,\ldots, \bfxi_\delta$ its reduced points. From \eqref{konk} and \eqref{tcr}, we have that
\begin{equation}\label{zizi}
\C[X_\sigma]/\langle F_1,\ldots, F_n\rangle\simeq \C[u_1,\ldots, u_n]^{G(\sigma)}/\langle \tilde{F}_1,\ldots, \tilde{F}_n\rangle\simeq \C[y_1,\ldots, y_\ell]/ I_{\sigma,\cA},
\end{equation}
and note that $\mbox{Res}_{F,\bfxi}$ is actually well defined in $\C[u_1,\ldots, u_n]^{G(\sigma)}/\langle \tilde{F}_1,\ldots, \tilde{F}_n\rangle,$ so it is an element of $(\C[y_1,\ldots, y_\ell]/ I_{\sigma,\cA})^\vee$ via \eqref{zizi}. 

We fix  $i\in\{1,\ldots, \delta\}.$ From Proposition \ref{cun}, we deduce that there exist $c_{i,1},\ldots, c_{i,\mu_i}\in\C$ such that
\begin{equation}\label{aunt}
\mbox{Res}_{F,\bfxi_i}=\sum_{j=1}^{\mu_i} c_{i,j}\partial_{i,j}.
\end{equation}
We show  below that the highest order coefficient cannot vanish.

\begin{proposition}\label{neq}
With notation as in \eqref{aunt}, we have that $c_{i,\mu_i}\neq0.$
\end{proposition}
\begin{proof}
Suppose that 
\begin{equation}\label{an}
\mbox{Res}_{F,\bfxi_i}=\sum_{j=1}^{\mu_i-1} c_{i,j}\partial_{i,j}.
\end{equation}
 As $\{\partial^{i,1},\ldots, \partial^{i,\mu_i}\}$ is a basis of $V_i,$ there exists $G\notin Q_i$ such that $\partial_{i,j}(G)=0$ if $j<\mu_i,$ and $\partial_{i,\mu_i}(G)\neq0.$   Via \eqref{zizi} we assume w.l.o.g. that $G\in\C[u_1,\ldots, u_n]^{G(\sigma)}.$ As $G\notin Q_i,$ from Proposition \ref{lresid} we deduce that  there is $H\in\C[u_1,\ldots, u_n]$ such that
\begin{equation}\label{deux}
\Res_{\bf0}\left(\frac{G \cdot H}{\tilde{F}_{1}\ldots \tilde{F}_n}\,d u_1\wedge\ldots\wedge d u_n\right)\neq0.
\end{equation}
As $H$ is not necessarily in $\C[u_1,\ldots, u_n]^{G(\sigma)},$ we will replace it with 
$$H^\sigma:=\sum_{g\in G(\sigma)} H(g(u))\, g(u_1)\ldots g(u_n)\in\C[u_1,\ldots, u_n]^{G(\sigma)},$$ and note that
$$\begin{array}{ccl}
\Res_{\bf0}\left(\frac{G \cdot H^\sigma}{\tilde{F}_{1}\ldots \tilde{F}_n}\frac{d u_1}{u_1}\wedge\ldots\wedge \frac{d u_n}{u_n}\right)& = &
\sum_{g\in G(\sigma)}\Res_{\bf0}\left(\frac{G \cdot H(g(u))}{\tilde{F}_{1}\ldots \tilde{F}_n}d g(u_1)\wedge\ldots\wedge dg( u_n)\right)\\
&=& |G(\sigma)|\Res_{\bf0}\left(\frac{G \cdot H}{\tilde{F}_{1}\ldots \tilde{F}_n}d u_1\wedge\ldots\wedge d u_n\right)\neq0,
\end{array}
$$
so we actually have that $\Res_{F,\bfxi_i}(G\,H^\sigma)\neq0.$ 

In addition, by construction we have $H^\sigma(\bf0)=0.$ 
This consideration, combined with Leibniz's rule \eqref{lbtz} applied to $\partial_{i,j}(G\cdot H^\sigma),$ and using that the subspace $V_i$ is closed,  and that $\partial_{i,j}(G)=0$ if $j<\mu_i,$  shows that 
 $\partial_{i,j}(G\cdot H^\sigma)=0$ for all $j<\mu_i.$ But this implies that 
 $$0\neq \Res_{F,\bfxi_i}(G\,H^\sigma)=\sum_{j=1}^{\mu_j-1}c_{i,j} \partial_{i, j}(G\,H^\sigma)=0,$$ 
a contradiction.
\end{proof}

\begin{corollary}\label{colau}
If $\{\partial_{i,1},\ldots, \partial_{i,\mu_i}\}$ is a consistently ordered basis of  $\Big(\C[y_1,\ldots, y_\ell]/ Q_i\Big)^\vee,$ then so is $\left\{\partial_{i,1},\ldots, \partial_{i,\mu_i-1},\Res_{F,\bfxi_i}\right\}.$
\end{corollary}

\section{The proofs of our main results}\label{sec:proofs}
In this section we state and prove Theorem~\ref{oa}. The input data is a system $F_1,\ldots, F_n\in S$ of homogeneous polynomials of degrees $\alpha_1,\ldots, \alpha_n,$ such that the lattice polytopes $P_{\alpha_1},\ldots, P_{\alpha_n}$ in ~\eqref{pd} define an indecomposable sequence, and such that  $V_X(F_1,\ldots, F_n)$   is zero-dimensional.
We then use Theorem~\ref{oa} to present the proofs of our main Theorems~\ref{smt} and~\ref{fonlyif} stated in the Introduction.

The proof of the following lemma is straightforward.

\begin{lemma}\label{tec}
Let $X= X_\Sigma$ be a complete toric variety and $X'$ be a complete simplicial toric variety associated to a simplicial refinement of $\Sigma$ without adding new rays. Call $\pi: X'\to X$ the associated proper map.
Then
\begin{itemize}
\item For any $x \in X$ in the simplicial part of $X$ (as in Definition~\ref{def:simplicial}),   $\pi^{-1}(x)$ consists of a single point.
\item 
Let $S$ be the common homogeneous coordinate ring of $X$ and $X'$ . Given
 a homogeneous system $F_1,\ldots, F_n\in S$ with common zeros  contained in the simplicial part of $X$,  the common zero locus $\{ x \in X \, : \, F_i(x)=0, i=1, \dots n \}$ is finite in $X$ if and only if $\{ x \in X' \, : \, F_i(x)=0, i=1, \dots n \}$ is finite in $X'$.
 \end{itemize} 
\end{lemma}

As each of the polytopes $P_{\alpha_1},\ldots, P_{\alpha_n}$  has all its vertices in $M,$  by Bernstein Theorem \cite{ber75}, the degree or $V_X(F_1,\ldots, F_n)$  is equal to  $w:=MV(P_{\alpha_1},\ldots, P_{\alpha_n}),$ the mixed volume of these polytopes. 
Set $I:=\langle F_1,\ldots, F_n\rangle,$ and $V_X(I)=\{\bfxi_1,\ldots, \bfxi_v\}.$ In addition, for each $i=1,\ldots, v,\,$ let $\mu_i$ be the multiplicity of $\bfxi_{i}$ as a zero of $V_X(F_1,\ldots, F_n).$ We have that $\mu_1+\ldots+\mu_v=w.$ 

Let $\alpha$ be any divisor, and consider the following morphism which appears in the proof of Corollary $3.1$ in \cite{BT22}:
\begin{equation}\label{square}
 (S/I)_{\alpha}  \stackrel{\psi_\alpha}{\longrightarrow}
\C^w 
\end{equation}
where the standard basis of  $\C^w$ is denoted with $e_{i,j},\ 1\leq i\leq v,\, 1\leq j\leq \mu_i,$  and $\big(\psi_{\alpha}([H])\big)_{(i,j)}=\partial_{i,j}(\tau_i(H^{\sigma_ i})),\ \partial_{i,j}$ being defined in \eqref{pij} with $\partial_{i, \mu_i} = \Res_{F,\bfxi_i}$ as in Corollary~\ref{colau}, the cone $\sigma_i$ verifies that $\bfxi_i\in X_{\sigma_i},\, H^{\sigma_i}$ is the dehomogenization of $H$ with respect to the cone $\sigma_i$ as in \eqref{deho}, and
\begin{equation}\label{taui}
\tau_i:\C[X_{\sigma_i}]\to\C[y_1,\ldots, y_{\ell_i}]/ I_{\cA_i}
\end{equation} is the ring isomorphism  \eqref{konk}.

The following two auxiliary results will be of use for the proof of Theorem \ref{oa} below.

\begin{lemma}\label{haux}
Let $T_1$ and $T_2$ matrices such that $T_1$ has maximal rank and $T_2$ has rank $1.$ Then, $T_1+T_2$ is either invertible or has corank $1.$
\end{lemma}

\begin{proof}
The subadditivity of the rank implies that $\mbox{rk}(T_1)\leq \mbox{rk}(T_1+T_2)+\mbox{rk}(-T_2)$, from which the result follows.
\end{proof}

\begin{lemma}\label{hhaux}
If $X$ is a complete normal toric variety, $\alpha_1,\ldots, \alpha_n\in A_{n-1}(X)$ are the degrees of nef Cartier divisors such that that the sequence $P_{\alpha_1},\ldots, P_{\alpha_n}$ is indecomposable,  and $F_1,\ldots, F_n\in S$ homogeneous  of respective degrees $\alpha_1,\ldots, \alpha_n$ with $V_X(F_1,\ldots, F_n)$ finite, then there is a degree $\alpha_0\in A_{n-1}(X)$ such that both $\alpha_0$ and $\alpha_0-\rho_0$ are the degrees of $\Q$-Cartier divisors, and moreover,
\begin{enumerate}
\item $P_{\alpha_0}$ is $n$-dimensional,
\item there exists $F_0\in S_{\alpha_0}$ which does not vanish on the points of $V_X(F_1,\ldots, F_n).$
\item there exists $F^*_0\in S_{\alpha_0-\rho_0}$ which does not vanish on the points of $V_X(F_1,\ldots, F_n).$
\end{enumerate}
\end{lemma}
\begin{proof}
Let $\beta_0\in A_{n-1}$ be such that it is the degree of a $\Q$-Cartier divisor, with $P_{\beta_0}$ being $n$-dimensional, 
and set $\alpha_0:=\beta_0+\rho_0.$ Then it is clear that both $\alpha_0$ and $\beta_0=\alpha_0-\rho_0$ are $\Q$-Cartier. The result then follows from Corollary $4.1$ in \cite{BT22}.
\end{proof}

\begin{theorem}\label{oa}
Let $X$ be a complete toric variety, and $\alpha_1,\ldots, \alpha_n\in A_{n-1}(X)$ degrees of nef Cartier divisors, such that the sequence of lattice polytopes $P_{\alpha_1},\ldots, P_{\alpha_n}$ is indecomposable. 
Consider $F_1,\ldots, F_n\in S$ homogeneous  of respective degrees $\alpha_1,\ldots, \alpha_n$ with $V_X(F_1,\ldots, F_n)$ 
finite and contained in the simplicial part of $X$. Then, the following assertions hold: 
\begin{itemize}
\item[(i)] Let  $\alpha=\alpha_1+\ldots+\alpha_n-\rho_0$.  The image of the map $\psi_\alpha$ in~\eqref{square} has codimension one in $\C^w$. It lies in the hyperplane defined by the sum of the coordinates $(i, \mu_i)$ with $i=1, \dots, v$, equal to $0$.
\item[(ii)] The hyperplane in item i) is the image  of all  $G\in S_\alpha$ such that
$$\sum_{i=1}^v\mbox{Res}_{\bfxi_i}\Big(\frac{G\,\Omega}{F_1\ldots F_n}\Big)=0.
$$

\item[(iii)]
Assume further that $\sigma\in\Sigma$ is a maximal simplicial cone such that  $\emptyset\neq V_{X_\sigma}(F_1^\sigma, \ldots, F_n^\sigma)$ and
$V_{X_\sigma}(F_1^\sigma, \ldots, F_n^\sigma)\subsetneq V_{X}(F_1,\ldots, F_n)$. 
Then, considering $\C[X_\sigma]$ as the $0$-th graded piece of the localization of $S$ at the monomial $\bfx^{\widehat{\sigma}}$, any $H\in\C[X_\sigma]$ has a representative  modulo $\langle F_1^\sigma,\ldots, F_n^\sigma\rangle\subset\C[X_\sigma]$ with support lying in the interior of the translated of $P$,
 where the vertex corresponding to $\sigma$ gets mapped to the origin in the translation.
 \end{itemize}
\end{theorem}

\begin{proof}
As $V_X(F_1, \dots, F_n)$ lies in the simplicial part of $X$, we can replace $X$ by $X'$ and- thanks to Lemma \ref{tec}- we can then assume w.l.o.g. that $X$ is a simplicial toric variety,  so the local residue
$\mbox{Res}_{\bfxi_i}\Big(\frac{G\,\Omega}{F_1\ldots F_n}\Big)$ at any common zero $\bfxi_i$  is defined thanks to the results of Section~\ref{sss}.
We know that if $\{\partial_{i,1},\ldots, \partial_{i,\mu_i}\}$ is 
 a consistently ordered basis of the dual of the local ring at $\bfxi_i$, where $\partial_{i,\mu_i}$ is equal to the local residue (as in Corollary \ref{colau}), then
$$\sum_{i=1}^v \big(\psi_{\alpha}([G])\big)_{(i,\mu_i)}=\sum_{i=1}^v\partial_{i,\mu_i}(\tau_i(G^{\sigma_ i}))=
\sum_{i=1}^v\mbox{Res}_{\bfxi_i}\Big(\frac{G\,\Omega}{F_1\ldots F_n}\Big)=0,$$
by Corollary \ref{ejsi}.  
We emphasize that the  fact that $\sum_{i=1}^v\mbox{Res}_{\bfxi_i}\Big(\frac{G\,\Omega}{F_1\ldots F_n}\Big)=0$ 
holds independently of whether the basis $(\partial_{i,j})$ is consistently ordered or not, and each of the local residues 
$\mbox{Res}_{\bfxi_i}\Big(\frac{G\,\Omega}{F_1\ldots F_n}\Big)$ is  a non trivial $\C$-linear combination of 
the $\partial_{i,j}(\tau_i(G^{\sigma_ i})), 1\leq j\leq \mu_i$ for any basis.  So, to prove items (i) and (iii) it is enough to show that the image of $\psi_\alpha$ has codimension one in $\C^w.$

To do so, pick $\alpha_0$ and $F_0$ as in Lemma \ref{hhaux}, and set  $\beta:=\alpha_0+\alpha_1+\ldots+\alpha_n-\rho_0$.  
From Lemma \ref{hhaux} we deduce that this value of $\beta$ belongs to the regularity of $I$ according to Theorem $4.4$ in \cite{BT22},  
and hence we have that $\dim_\C(S/I)_{\beta}=w.$ Moreover, from Corollary $3.1$ in loc. cit, we also deduce that $\psi_\beta$ defined in \eqref{square} is an isomorphism. 
So, for each pair $(i,j),\, 1\leq i\leq v,\, 1\leq j\leq \mu_i,$ there is an element $B_{i,j}\in S_\beta$ such that its class $[B_{i,j}]\in S/I$ gets mapped to $e_{i,j}\in\C^w.$

From the hypothesis, as $P_{\alpha_0}$ is $n$-dimensional and the sequence $P_{\alpha_1},\ldots, P_{\alpha_n}$ is indecomposable, 
we deduce from Proposition \ref{cd1} that $\dim_\C\left(S/\langle F_0,\ldots, F_n\rangle\right)_{\rho_\bfF}=1.$ Let $\Delta\in S_{\rho_\bfF}$ 
be an element whose class is not zero in this quotient. We get then that
$[B_{i,j}]=[A_{i,j}F_0+\lambda_{i,j}\Delta]$ for suitable $A_{i,j}\in S_\alpha,$ and $\lambda_{i,j}\in\C.$ 

So we have that, for any pair $(k,\ell)$, by considering the ring isomorphism $\tau_k$~\eqref{taui},
\begin{equation}\label{ozu}
\big(\psi_\beta([B_{i,j}])\big)_{(k,\ell)}=\partial_{k,\ell}\big(\tau_k(A^{\sigma_k}_{i,j}F^{\sigma_k}_0+\lambda_{i,j}\Delta^{\sigma_k})
\big)=
\left\{\begin{array}{ccc}
1 & \mbox{if} & (k,\ell) = (i,j),\\
0& \mbox{if} & \mbox{not.}
\end{array}
\right.
\end{equation}
Note that  $\alpha$ is not necessarily the degree of a Cartier divisor, but $\alpha+\alpha_0$ is Cartier, so $A^{\sigma_k}_{i,j}$ 
and $\Delta^{\sigma_k}$ should be defined as follows: if $P_{\alpha+\alpha_0}=\{m\in M_\R:\,\langle m,\eta_j\rangle\geq -a_j,\ j=
1,\ldots, n+r\},$ with vertices $m_\sigma$ for maximal $\sigma\in\Sigma,$ then 
$$A^{\sigma_k}_{i,j}:=\frac{x_1\ldots x_{n+r}A_{i,j}}{\prod_{p=1}^{n+r}x_p^{\langle m_{\sigma_k},\eta_p\rangle+a_p}}, 
\ \Delta^{\sigma_k}:=\frac{x_1\ldots x_{n+r}\Delta}{\prod_{j=p}^{n+r}x_j^{\langle m_{\sigma_k},\eta_p\rangle+a_p}}.$$
From Proposition \ref{goxi}, we deduce that $\Delta\in \langle F^{\sigma_k}_0,\ldots, F^{\sigma_k}_n\rangle$ and hence, we have $\Delta^{\sigma_k}=B_k F^{\sigma_k}_0,$ 
for a suitable $B_{k}\in S_{\bfx^{\widehat{\sigma_k}}}$,  modulo  
$\langle F^{\sigma_k}_1,\ldots, F^{\sigma_k}_n\rangle$.
So, we get
\begin{equation}\label{sale}
\partial_{k,\ell}\big(\tau_k(A^{\sigma_k}_{i,j}F^{\sigma_k}_0+ \lambda_{i,j}\Delta^{\sigma_k})
\big)=
\partial_{k,\ell}\big(\tau_k(A^{\sigma_k}_{i,j}+\lambda_{i,j}B_k)\tau_k(F^{\sigma_k}_{0})\big).
\end{equation}
Let $M_A, T_0, D_0$ and $C$ be the $w \times w$ matrices whose rows and columns are indexed by the pairs $(i,j),\, 1\leq i\leq v,\, 1\leq j\leq \mu_i$, defined as follows:
\begin{itemize}
\item 
$M_A:=\left(\psi_{\alpha}([A_{i,j}])\right)_{(k,\ell)}=\left(\partial_{k,\ell}\big(\tau_k(A^{\sigma_k}_{i,j}))\right)_{(i,j),\,(k,\ell)}$;

\item $T_0$ has a diagonal-block structure with $v$ blocks, the one indexed by $k\in\{1,\ldots, v\}$ being the 
triangular $w_k \times w_k$-matrix appearing in \eqref{triang} with $\tau_k(F^{\sigma_k}_0)$ instead of $g,$ and the rest of
 its structure being given by  the behaviour of the multiplication by $\tau_k(F^{\sigma_k}_0)$ in $\C[X_{\sigma_k}]$ around the point $\bfxi_k;$
\item $D_0$ is diagonal, having the $\lambda_{i,j}$'s in the diagonal;
\item $C$  has all the elements in the column indexed by $(k,\ell)$ equal to the constant value $\partial_{k,\ell}(\tau_k(B_k)).$ 
\end{itemize}
As $F_0$ does not vanish on $V_X(F_1,\ldots, F_n),$  we deduce that $T_0$ is triangular and of maximal rank. Moreover, 
from \eqref{ozu} and \eqref{triang}, we deduce that
$$T_0\cdot(M_A+D_0\cdot C)=I,
$$
with $I$ the identity matrix.  By inverting $T_0$, we see that $M_A$ is equal to a matrix of the form
 $T_1+T_2$ as in Lemma \ref{haux}. So, we deduce that $M_A$ has either maximal rank or corank $1.$ 
As $M_A$ factorices through $\psi_\alpha$ which we know cannot have maximal rank, we deduce then that $M_A$ 
has corank $1$, and a fortiori that the image of $\psi_\alpha$ also has corank $1$. This proves item~(i).

To prove item~(iii), assume w.l.o.g. that $V_{X}(F_1,\ldots, F_n)\cap X_\sigma=\{\bfxi_1,\ldots,\bfxi_\ell\}$ and that $\bfxi_{\ell+1}\notin X_\sigma.$ Set  $w_\sigma:=\mu_1+\ldots+\mu_\ell.$ From \eqref{tcr} we have that $w_\sigma$ is equal to the 
$\C$-dimension of $\C[X_\sigma]/\langle F_1^\sigma,\ldots, F_n^\sigma\rangle$.  Note that the restriction of  $\psi_{\alpha}$ to $\C^{w_\sigma}$ is onto as the operator
$\mbox{Res}_{\bfxi_{\ell+1}}\Big(\frac{--\,\Omega}{F_1\ldots F_n}\Big)$ can be 
expressed as a nontrivial linear combination of the duals of $e_{\ell+1, j},\, 1\leq j\leq \mu_{\ell+1},$ i.e. they do not depend 
on the first $w_\sigma$ vectors of the standard basis of $\C^w.$
This implies that there are homogeneous forms $H_1,\ldots, H_{w_\ell}\in S_{\alpha-\rho_0}$ which are mapped to $e_{1,1},\ldots, e_{\ell,\mu_\ell}\in\C^{w_\sigma}.$ Dehomogenizing $H_j\,x_1\ldots x_{n+r},\,1\leq j\leq w_\ell$ to $\C[X_\sigma]$ via \eqref{deho} we obtain a maximal linearly independent set in the quotient ring 
$\C[X_\sigma]/\langle F_1^\sigma,\ldots, F_n^\sigma\rangle$ made by polynomials with support as in the claim. 
\end{proof}

\begin{remark}
As a straightforward consequence of Theorem \ref{oa}, we deduce that $\dim(S/I)_\alpha\geq w-1$
 for  $\alpha=\alpha_1+\ldots+\alpha_n-\rho_0,$ as it follows from
Proposition \ref{bound} for the codimension one case.
\end{remark}

\subsection{Proofs of the main theorems}

We are now ready to prove Theorem~\ref{smt}, that we restate for convenience here. Let $f_1,\ldots, f_n\in\C[\bft^{\pm1}]$ be as in \eqref{eq:efes} with  $P_1,\ldots, P_n$ indecomposable. Assume that $V_T(\bff)$ is nonempty and that the intersection of the closures of the hypersurfaces $(f_i=0)$ for $i=1,\dots,n$, has dimension zero in  the complete toric variety $X_P$ associated with $P=P_1+\dots +P_n$ and it is contained in the simplicial part of $X$.
Then, the following are equivalent:
\begin{enumerate}
\item[i)] $\deg(V_T(\bff))=\MV(P_1,\ldots, P_n)$.
\item[ii)] For any $ h_0\in\C[\bft^{\pm1}]$ with support contained in $P^\circ$, 
$\mbox{Res}^T_{\bff}(h_0) =0$.
\item[iii)] There is no Laurent polynomial  $p_J\in  \C[\bft^{\pm1}]$ with support contained in $P^\circ$ such that $J^T_\bff\equiv p_J$ modulo the ideal 
$\langle f_1,\ldots, f_n\rangle$ in $\C[\bft^{\pm1}]$. 
\end{enumerate}

\begin{proof}[Proof of Theorem \ref{smt}] $^{}$
If item i) is satisfied, then the validity of item ii) is the Euler-Jacobi Theorem~\ref{koo} (proven in~\cite[Corollary $5$]{CD97}).
Also item ii) easily implies item iii), because $\mbox{Res}^T_{\bff}(J^T_\bff) \neq 0$ since we are assuming that $V_T({\bf f}) \neq \emptyset$ and then $\mbox{Res}^T_{\bff}(J^T_\bff)$ equals a sum of positive integers (the multiplicities at all the zeros in the torus).

To prove the remaining implication iii) $\implies$ i), we show now that
 if item i) does not hold then item iii) does not hold. Denote by $F_1, \dots, F_n$ are the homogenizations of the input polynomials. By hypothesis, ${\rm dim}(V_{X_P}(F_1, \dots, F_n) )=0$.
 Note that the toric jacobian $J^T_\bff$ is  supported in $P\cap\Z^n,$ so it can be $P$-homogenized. 
Passing through a refinement of cones without adding new rays, we may assume w.l.o.g. that $X_P$ is simplicial (Proposition~11.1.7 in \cite{CLS11}) keeping our hypothesis that
there is a finite number of zeros by Lemma~\ref{tec}. Recall also that we would not be changing the homogeneous coordinate ring. 
We are assuming that
$0 < \deg(V_T(\bff)) <\MV(P_1,\ldots, P_n)$.
 Let ${\bf \xi} \in X_P \setminus T$. By Proposition~\ref{prop:1}, we have that $T = \cap_\sigma (X_P)_\sigma$. Then, there exists a (simplicial) cone $\sigma$ such that ${\bf \xi} \notin {X_P}_\sigma$. Then,
 $\emptyset\neq V_{{X_P}_\sigma}(F_1^\sigma, \ldots, F_n^\sigma)\subsetneq V_{X_P}(F_1,\ldots, F_n)$. 
We apply the third item in Theorem~\ref{oa} to conclude that the homogenization of $J^T_\bff$ can be represented modulo $\langle F_1^\sigma,\ldots, F_n^\sigma\rangle$ 
as a polynomial with exponents lying in the interior of $P$. Going back to the coordinates in the torus via \eqref{tx}, we have that item iii) is not valid.
\end{proof}

We now present the proof of our second main theorem stated in the Introduction.

\begin{proof}[Proof of Theorem~\ref{fonlyif}]
From Proposition \ref{prop:amplefull}, we deduce that $P_{\alpha_i},\, i=0,\ldots, n,$ are lattice polytopes, and moreover that this family is indecomposable, as each of these polytopes is full dimensional. Assume w.l.o.g. that $i=0.$ Set as before  $I=\langle F_1,\ldots, F_n\rangle,$ and
$V_X(F_1,\ldots, F_n)=\{\bfxi_1,\ldots, \bfxi_v\},$ contained in the simplicial part of $X$ and having degree $w\geq v.$
Suppose also w.l.o.g. that $F_0(\bfxi_1)=0.$   

Let $X'$ be a simplicial desingularization of $X$ without the addition of rays. Note that the definition of $J_\bfF$ does not depend on this refinement but the definition of $\Delta_{\bfF,\bfsigma}$
 is based on the original fan $\Sigma,$ so some care must be taken with this case.

As the divisor $\rho_\bfF$
fits within the hypothesis of Theorem 1.2 in \cite{BT22}, we have that $\dim_\C(S/I)_{\rho_\bfF}=w,$ and  $\psi_{\rho_\bfF}$ is an isomorphism. Hence, 
$\psi_{\rho_\bfF}([\Delta_{\bfF,\sigma}])$ completely determines $\Delta_{\bfF,\sigma}$ modulo $I$.
We will pick a dual basis $(\partial_{ij})_{(i,j)}$ such that  that for each $i\in\{1,\ldots, v\},$
 $\{\partial_{i,1},\ldots, \partial_{i,\mu_i}\}$ is a consistently ordered basis  
 with $\partial_{i,\mu_i}=R_{\bff_\sigma, \bfxi_i}\Big[\frac{d\bfy}{f_{1,\sigma}\ldots f_{\ell,\sigma}}\Big] ,$ as in Corollary \ref{colau}.

Fix $k$ and recall that we are assuming that the cone $\sigma_k$ is simplicial. From Proposition \ref{oxi}, 
we deduce that $\Delta^{\sigma_k}_{\bfF,\sigma}= B_{k}\cdot F^{\sigma_k}_{0}$ modulo $\langle F^{\sigma_k}_1,\ldots, F^{\sigma_k}_n\rangle,$
for a suitable $B_k\in\C[X_{\sigma_k}].$ Then, as in \eqref{sale}, we have  for any index $\ell$ with ${1\leq\ell\leq\mu_k}$:
\begin{equation}\label{oxig}
\partial_{k,\ell}\big(\tau_k(\Delta^{\sigma_k}_{\bfF,\sigma})
\big)=
\partial_{k,\ell}\big(\tau_k(B_{k})\cdot u_k(F^{\sigma_k}_{0})\big)
=M_{\tau_k(F^{\sigma_k}_{0})}\cdot \big(\partial_{k,\ell}(\tau_k(B_{k}))\big),
\end{equation}
where $M_{\tau_k(F^{\sigma_k}_{0})}$ is the triangular matrix appearing in \eqref{triang} 
with $g$ replaced by $u_k(F^{\sigma_k}_0)$ in the columns indexed by $(k,\ell),\,1\leq\ell\leq \mu_k.$ 
We deduce then that $\psi_{\rho_\bfF}([\Delta_{\bfF,\sigma}])$ is a linear combination of the columns of $M_{\tau_k(F^{ \sigma_k}_0)}, \ k=1,\ldots, v.$

From Theorem \ref{oa}, and due to the fact that we are working with a consistently ordered basis 
having the local residue as their last elements,  we deduce that it is possible to find  for each $(i,j)$ with $j<\mu_i,$  a form  $B_{i,j}\in S_{\alpha_1+\ldots+\alpha_n-\rho_0}$ such that
$\psi_{\alpha_1+\ldots+\alpha_n-\rho_0}([B_{i,j}])=e_{i,j};$ and also for each $i=2,\ldots, v,$ a form  
$B_{i,\mu_i}\in S_{\alpha_1+\ldots+\alpha_n-\rho_0}$ such that
$\psi_{\alpha_1+\ldots+\alpha_n-\rho_0}([B_{i,\mu_i}])=e_{i,\mu_i}-e_{1,\mu_1},$ as in the dual basis of the standard of $\C^w,$ we have that the sum of the coordinates $(i, \mu_i)$ for $i=1, \dots,v$ is equal to $0$.

Then, for all $(i,j)\neq(1,\mu_1)$:
\begin{equation}\label{avv}
\psi_{\rho_\bfF}([B_{i,j}F_0])=M_{F_0}\cdot \left(\psi_{\alpha_1+\ldots+\alpha_n-\rho_0}([B_{i,j}])\right),
\end{equation}
which is either the column $(i,j)$ of $M_{F_0}$ (if $j<\mu_i$), or the column $(i,\mu_i)$ minus the column $(1,\mu_1)$ of this matrix. 
But  as $F_0(\bfxi_1)=0,$ this  column is identically zero. So, we deduce from above that the columns 
of~\eqref{avv} generate the $\C$-vector space of the columns of $M_{F_0}.$
As \eqref{oxig} is a linear combination of these columns, we then deduce that there exist $d_{i,j}\in\C$ such that
$$\psi_{\rho_\bfF}\left(\Big[\Delta_{\bfF,\sigma}-\Big(\sum_{(i,j)\neq(1,\mu_1)} d_{i,j}B_{i,j}\Big)F_0\Big]\right)=0,
$$
which proves the claim for $\Delta_{\bfF,\sigma}.$ A similar result as in \eqref{oxig} holds if we replace this polynomial with 
$J_{\bfF}$ in the cases where all the degrees are a common multiple of a fixed one. The result then follows again mutatis mutandis the proof above.
\end{proof}

\bigskip
\section{Open questions}\label{sec:open}

We end the article with some questions from further work. 
The first general question is the following obvious one:
\begin{question}\label{q:one}
Is it possible to get rid of the hypothesis that we used all along this paper on the finiteness of $V_X(F_1,\ldots, F_n)$ (and its containment in the simplicial part of $X$)?
\end{question}

The answer to Question~\ref{q:one} is true in the case of projective space. Assume that $f_1, \dots, f_n$ are polynomials of degrees $d_1, \dots, d_n$ with finitely many zeros in the torus with the sum of their multiplicities $v$ being positive but less than the B\'ezout number $d_1\dots d_n$. If there is an infinite number of zeros in $\P^n$, then there are infinitely many common solutions
to their homogenizations $F_1, \dots, F_n$ at the union of the divisors  $(x_i=0)$ for $i=0, \dots, n$. We assume without loss of generality that this happens at $(x_0=0)$. But then,  the homogeneous polynomials $G(f_1),\dots, G(f_n)$ of  degree $d_1, \dots, d_n$ in $f_1, \dots, f_n$ respectively, don't have isolated zeros in $\C^n$. W. Vasconcelos proved in~\cite{vas92} that the colon ideal of  polynomials sending the jacobian $J_{G(f)}$ into the ideal $I_G = \langle G(f_1), \dots, G(f_n)\rangle$  equals the
intersection of the primary components of the $I_G$  corresponding to its
isolated primary components of dimension $0$. 
It follows that in this case, the colon ideal equals the whole polynomial ring $\C[x_1, \dots, x_n]$, and so 
$J_{G(f)} \in I_G$. It is easy to see that then the standard jacobian $J_f$ is equivalent
modulo the ideal $\langle f_1, \dots, f_n \rangle$ to a polynomial of degree smaller than $\sum_{i=1}^n d_i - n$. But then the global residue $\mbox{Res}^T_\bff(J_f) =v \not= 0$, proving the converse to the toric Euler-Jacobi Theorem~\ref{koo}.

\smallskip

We also pose the following questions:

\begin{question}
Let  $F_1, \dots,F_n$ and $X$ be as in the statement of Theorem \ref{oa}. We deduce that if the zeroes of the ideal  $I$  generated by $F_1, \dots,F_n$ are finite in $X$, then  the quotient  $(S/I)_{\alpha}$ has dimension at least $MV(P_{\alpha_1},\ldots, P_{\alpha_n})-1$ for $\alpha=\alpha_1+\ldots+\alpha_n-\rho_0$.  By standard properties of Hilbert functions, we know that the dimension of this quotient is at most $MV(P_{\alpha_1},\ldots, P_{\alpha_n}),$ and in all the classical cases, it is known that the dimension is actually the lower bound. Does this hold in general?
\end{question}

\begin{question}\label{question3}
Let $f_1, \dots, f_n$ as in~\eqref{eq:efes} with respective supports in a sequence of indecomposable polynomials 
$P_1, \dots, P_n$ with $P=P_1 + \dots + P_n$. Are the following statements equivalent?:
\begin{enumerate}
\item $V_{T}(f_1,\ldots, f_n)$ is finite of degree strictly less than $MV(P_1,\ldots, P_n)$.
\item The quotient ring $\C[\bft^{\pm1}]/\langle f_1,\ldots, f_n\rangle$ can be generated by the monomials lying in the interior of $P$.
\end{enumerate}
\end{question}

Note that $(2)\implies (1)$ in Question~\ref{question3} follows from the Toric Euler-Jacobi Theorem (Theorem \ref{koo}). 
In addition, if  the intersection of the closures of the hypersurfaces $(f_i=0)$ for $i=1,\dots,n$, has dimension zero in $X_P$,  the proof of Theorem~\ref{smt}, shows that $(1)$ implies that any Laurent polynomial $h$ with support contained in $P$, has a representative in the quotient in (2) in terms of monomials in the interior of $P$ (and not only the toric jacobian). Moreover, this is true if the support of $h$ is contained in one of the simplicial cones of $\Sigma$.

We also have that $(1)\implies (2)$ if $\langle f_1,\ldots, f_n\rangle$ is radical in $\C[\bft^{\pm1}]$  and there is a simplicial cone $\sigma$ such that $V_{\T}(f_1,\ldots, f_n)\subset V_{X_\sigma}(F^\sigma_1,\ldots, F^\sigma_n)\subsetneq V_{X}(F_1,\ldots, F_n),$ with no positive dimensional components in $V_{X}(F_1,\ldots, F_n).$ Indeed, for any $\bfxi\in V_{\T}(f_1,\ldots, f_n)$, there is a polynomial $g_\bfxi\in\C[X_\sigma]$ supported in the interior of $P$ whose value in $\bfxi$ is equal to $1,$ and when evaluated in the rest of the elements in the finite set $V_{\T}(f_1,\ldots, f_n)$ is equal to zero. 
 This is because  $\bfxi$ is also a single root of $V_{X_\sigma}(F^\sigma_1,\ldots, F^\sigma_n),$ and hence the set of evaluations on all $\bfxi$ in this set is linearly independent in $\left(\C[X_\sigma]/\langle F^\sigma_1,\ldots, F^\sigma_n\rangle\right)^\vee.$ One can extend this set to a basis of this space, and then compute its dual basis in $\C[X_\sigma]/\langle F^\sigma_1,\ldots, F^\sigma_n\rangle,$ the latter being generated by monomials lying in the interior of a translate of $P$ thanks to Theorem \ref{oa} iii). Denote with $\{g_\bfxi,\,\bfxi\in V_\T(f_1,\ldots, f_n)\}$ the first elements of such a basis. Note that the support of all of these polynomials is also contained in the interior of a translate of $P$.  By Hilbert's Nullstellensatz, we deduce that this set generates  $\C[\bft^{\pm1}]/\langle  f_1,\ldots, f_n\rangle,$ as for any $g\in\C[\bft^{\pm1}]$ there is a linear combination of the $g_\bfxi$ which has the same value of $g$ in  $V_{\T}(f_1,\ldots, f_n),$ and $\langle f_1,\ldots, f_n\rangle$ is radical.

\end{document}